\documentclass[preprint,12pt,authoryear]{elsarticle}

\usepackage{amsmath}
\usepackage{amssymb}
\usepackage{hyperref}
\usepackage{enumitem}
\usepackage{graphicx}
\usepackage{subcaption}
\usepackage{booktabs}
\usepackage{tabularx}
\usepackage{makecell}
\usepackage{rotating}
\usepackage{multirow}
\usepackage{siunitx}
\usepackage{orcidlink}

\makeatletter
\def\ps@pprintTitle{%
 \let\@oddhead\@empty
 \let\@evenhead\@empty
 \def\@oddfoot{\footnotesize\itshape
  \textcopyright\ \the\year. Preprint submitted to arXiv \hfill \today}
 \let\@evenfoot\@oddfoot}
\makeatother

\begin{document}

\begin{frontmatter}

\title{An Efficient Continuous-Time MILP for Integrated Aircraft Hangar Scheduling and Layout}

\author[1]{Shayan Farhang Pazhooh\orcidlink{0009-0001-1366-7942}}
\ead{shayanfp2006mail@gmail.com}

\author[2]{Hossein Shams Shemirani\orcidlink{0000-0002-5255-7746}\texorpdfstring{\corref{cor1}}{}}
\ead{h.shams@iut.ac.ir}

\cortext[cor1]{Corresponding Author}

\affiliation[1]{%
  organization={Department of Industrial Engineering},
  addressline={Isfahan University of Technology},
  city={Isfahan},
  country={Iran}
}

\affiliation[2]{%
  organization={Industrial Engineering Group, Golpayegan College of Engineering, Isfahan University of Technology},
  addressline={Golpayegan 87717-67498},
  country={Iran}
}

\begin{abstract}
Efficient management of aircraft MRO hangars requires the integration of spatial layout with time-continuous scheduling to minimize operational costs. We propose a continuous-time mixed-integer linear program that jointly optimizes aircraft placement and timing, overcoming the scalability limits of prior formulations. A comprehensive study benchmarks the model against a constructive heuristic, probes large-scale performance, and quantifies its sensitivity to temporal congestion. The model achieves orders-of-magnitude speedups on benchmarks from the literature, solving a long-standing congested instance in 0.11 seconds, and finds proven optimal solutions for instances with up to 40 aircraft. Within a one-hour limit for large-scale problems, the model finds solutions with small optimality gaps for instances up to 80 aircraft and provides strong bounds for problems with up to 160 aircraft. Optimized plans consistently increase hangar throughput (e.g., +33\% serviced aircraft vs. a heuristic on instance RND-N030-I03), leading to lower delay penalties and higher asset utilization. These findings establish that exact optimization has become computationally viable for large-scale hangar planning, providing a validated tool that balances solution quality and computation time for strategic and operational decisions.
\end{abstract}

\begin{keyword}
Aircraft Maintenance \sep Hangar Scheduling \sep Spatial-Temporal Optimization \sep Mixed-Integer Linear Programming \sep Heuristic Algorithm \sep Decision Support Systems
\end{keyword}

\end{frontmatter}


\section{Introduction}
The global aircraft maintenance, repair, and overhaul (MRO) market is a multi-billion dollar industry in which operational efficiency directly translates to profitability. Within this complex ecosystem, the maintenance hangar emerges as a pivotal yet capacity-constrained asset, often serving as a critical bottleneck in the MRO value chain. Consequently, optimizing hangar space and time utilization presents a paramount challenge for airlines and service providers, with direct implications for fleet readiness and financial performance.

The hangar scheduling problem represents a complex spatial-temporal optimization challenge, falling within the broader class of Operations Routing and Scheduling Problems (ORSP) that integrate logistical movements with operational tasks \citep{jiang2025operations}. It involves two intertwined subproblems: a temporal scheduling problem (when to bring aircraft in and out) and a spatial allocation problem (where to park them). Decisions made in one domain directly constrain the other. For instance, the physical placement of one aircraft can block the movement of another, dictating the sequence in which they must depart. Similarly, the scheduled arrival and departure times determine which aircraft must share the hangar space simultaneously, thus constraining the possible spatial layouts.

While this complex problem has been addressed in the literature, existing exact methods often rely on computationally intensive discrete-time or event-based formulations. A prominent example is the work of \citet{qin2019}, who developed an efficient event-based model that performed well in many scenarios. However, its performance proved highly sensitive to temporal congestion. In one particularly challenging benchmark instance—featuring only nine aircraft but highly congested arrival times—the model failed to find an optimal solution within a one-hour time limit; a challenge that, as we will demonstrate, our proposed formulation overcomes in a fraction of a second.

This computational barrier has naturally led many researchers to pivot towards heuristic and metaheuristic approaches to find high-quality solutions in a reasonable timeframe \citep{zhou2023membership}. While valuable, such methods by design do not guarantee optimality. This has fostered a prevailing assumption in the field: that finding provably optimal solutions for large and congested hangar scheduling problems is computationally prohibitive.

This paper directly challenges that assumption by introducing a novel and highly efficient continuous-time MILP formulation. The main contributions of this work are threefold:

\begin{enumerate}
    \item \textbf{A Novel Continuous-Time Formulation:} We propose an efficient MILP model that treats time as a continuous variable. This approach fundamentally reduces the model's complexity compared to discrete-time and event-based counterparts, enabling the efficient solution of larger and more complex problem instances.
    
    \item \textbf{An Extensive and Multi-Faceted Computational Study:} We validate our model through a rigorous computational campaign that includes direct benchmarking against the state-of-the-art on published instances, large-scale scalability tests on problems with up to 160 aircraft, and a controlled sensitivity analysis on temporal congestion.
    
    \item \textbf{An Insight-Driven Decision Support Framework:} We present a custom-built visualization tool that transforms the model's numerical output into an interactive dashboard. This framework serves as a powerful analytical instrument, enabling managers to bridge the gap between advanced analytics and actionable operational execution.
\end{enumerate}

Collectively, this research provides a holistic and scalable framework for optimizing hangar operations. By proving that exact, optimal solutions are now computationally within reach for large-scale, real-world problems, we equip MRO planners with a more powerful tool for enhancing throughput, reducing costs, and gaining significant managerial insights.

\section{Literature Review}
The optimization of maintenance activities is a well-established field in operational research, with a vast body of literature dedicated to developing models that enhance efficiency and reliability across various industries \citep{de2020review, arts2025fifty}. As categorized by a recent review from \citet{dinis2025maintenance}, maintenance management faces three core challenges: capacity planning, spare parts management, and scheduling. Within this broad domain, aircraft maintenance presents a unique set of challenges due to the high value of assets, stringent safety regulations, and the complex interplay of operational constraints. Research in aircraft maintenance spans multiple decision levels, from strategic fleet-wide planning—such as maximizing fleet availability \citep{gavranis2015fleet}, maintenance routing \citep{sriram2003, qin2024service}, and long-term check scheduling \citep{andrade2021}—to tactical resource allocation, including personnel rostering \citep{belien2012} and technician assignment \citep{chen2017}. Our work is situated at the operational core of MRO activities: the aircraft hangar scheduling problem, which lies at the intersection of scheduling, layout planning, and packing problems, directly addressing the challenge of finite hangar capacity. This paper addresses the operational-level decision of assigning precise roll-in/out times and physical parking positions for incoming aircraft, considering spatial and temporal constraints.

The development of exact models for this problem has evolved from static, single-period spatial allocation to dynamic, multi-period scheduling. An early key work by \citet{qin2018parking} addressed the static aircraft parking stand allocation problem, focusing on maximizing profits and safety margins for a single day by selecting a subset of aircraft. This was extended into a comprehensive, multi-period model by \citet{qin2019}, which presents an integrated mathematical model for the hangar scheduling and layout problem. Their work is comprehensive, addressing both temporal and spatial dimensions with detailed considerations for aircraft geometry and movement blocking. However, their formulation is based on a discrete-time framework, where decision variables are indexed by specific time points. While this method is conceptually straightforward, it suffers from significant scalability issues. As the authors note, their model struggles to find optimal solutions for instances with as few as 9 aircraft within a one-hour time limit, necessitating a rolling horizon approach for larger, more realistic instances, which sacrifices global optimality.

This computational barrier has steered many researchers towards developing sophisticated heuristic and metaheuristic algorithms. These approaches trade the guarantee of optimality for computational speed, aiming to produce near-optimal solutions for larger or more complex instances. For example, \citet{zheng2020hybrid} modeled the integrated aircraft scheduling and parking problem as a generalization of the two-dimensional strip-packing problem (2D-SPP), proposing a hybrid simulated annealing and variable neighborhood search (HSARVNS) capable of solving instances with up to 250 aircraft. Similarly, \citet{zhou2023membership} tackled a variant of the problem with membership-based priorities, developing a custom heuristic based on Bacterial Foraging Optimization (BFO) for smaller instances. Even studies that begin with exact MILP formulations often pivot to heuristics for practical application; \citet{agostinho2025planning}, for instance, developed a detailed MILP for maintenance task planning but ultimately proposed a matheuristic to solve large-scale instances sequentially. Other works have focused on specific operational constraints, such as the critical issue of aircraft blocking during multi-stage movements, which \citet{chen2024stand} addressed using a discrete-time network flow model. While these heuristic and specialized models offer powerful solutions, they inherently sacrifice global optimality or a fully integrated perspective.

The core challenge of the hangar problem—managing high-value assets through a capacitated, shared transfer system—is not unique to aviation. A strong parallel exists in the maritime industry with syncrolift dry dock scheduling, where ships are moved between the sea and onshore service bays via a shared lift and railway system. \citet{guan2025syncrolift} address this by developing an MILP model to minimize ship waiting times, explicitly modeling the transfer system as a series of capacity-constrained segments. Their work underscores the critical role of the transfer system as a bottleneck and highlights the value of integrated scheduling, reinforcing the general importance of the problem structure we address. To overcome the scalability limitations of exact models, other approaches have emphasized different facets of the hangar scheduling problem. Recognizing the problem's similarity to two-dimensional irregular packing, \citet{niu2024two} developed a two-stage metaheuristic framework to maximize hangar utilization. Others have viewed it as a variant of the bin packing problem, as seen in the work of \citet{witteman2021}, who modeled maintenance task allocation.

The trend towards more realistic and integrated models increasingly situates specific operational problems within the broader airport ecosystem. For instance, significant research has focused on optimizing other critical airport bottlenecks, such as runway operations through collaborative flight scheduling \citep{jiang2024collaborative} and ground-side efficiency via the scheduling of autonomous electric vehicles for baggage transport \citep{zhang2025adaptive}. While these studies are vital for overall airport performance, our work addresses the maintenance hangar—a vital and distinct component of airport logistics with its own unique set of spatial-temporal challenges.

Separately, other models have advanced the state-of-the-art by incorporating greater realism or novel methodologies. Research has expanded to include integrated problems like technician assignment \citep{qin2020two} and scheduling for moving assembly lines, where concepts like project splitting are introduced to manage resource constraints and time windows \citep{lu2019resource}. A significant trend is the focus on managing uncertainty. To this end, recent works have employed artificial intelligence and computer science paradigms. For example, \citet{hu2021reinforcement} developed a reinforcement learning-driven strategy for long-term maintenance decisions under uncertainty. More recently, \citet{yang2025out} proposed an operating system-inspired framework using out-of-order execution to handle real-time disruptions across the entire planning, scheduling, and execution (PSE) pipeline. While these models are invaluable for capturing operational dynamics and uncertainty, they often rely on simulations or heuristic policies, and do not provide provably optimal plans for the core deterministic scheduling problem that underpins stable operations.

This review of the literature culminates in a clear and compelling research gap: the absence of an exact and scalable formulation for the integrated hangar scheduling problem. While existing models offer valuable insights, they consistently face a computational barrier, particularly under high temporal congestion, forcing a reliance on heuristics that inherently sacrifice global optimality. This paper directly confronts this challenge by proposing a novel continuous-time MILP model—a fundamental paradigm shift from discrete-time or event-based frameworks. By significantly reducing model complexity, this approach is designed to deliver both optimality and scalability. Our model is further distinguished by its comprehensive cost structure and its ability to handle initial hangar states. By uniting a robust exact method with a practical visualization tool, our work offers a powerful, integrated framework for optimizing real-world MRO operations.

\section{Problem Definition and Mathematical Model}
The problem is formulated as a mixed-integer linear program (MILP). The following sections detail the sets, parameters, and variables that define the model, followed by an in-depth explanation of the objective function and constraints.

\subsection{Problem Definition and Assumptions}
The problem addressed in this paper is the integrated scheduling and spatial allocation of aircraft within a maintenance hangar. The objective is to determine which new aircraft to accept, when they should arrive and depart, and where they should be parked to minimize total operational costs, which include penalties for rejections and delays.

The operational environment is defined by the following key assumptions, which form the basis for our mathematical model:

\begin{itemize}
    \item \textbf{Hangar Layout and Access:} We model the hangar as a rectangular area with a single, primary entrance/exit located along the edge corresponding to the highest Y-coordinates. This "open-front" configuration, which is consistent with modern, large-span hangar designs that favor a wide layout for structural integrity and operational efficiency \citep{liu2025}, dictates that all aircraft movements are sequential and primarily oriented along the Y-axis.

    \item \textbf{Movement Exclusivity:} To prevent logistical conflicts and ensure safety, no two movement events (i.e., an aircraft rolling in or rolling out) can occur simultaneously. A minimum time gap, denoted by $\varepsilon_t$, must be maintained between any two consecutive movements. This is enforced through a set of temporal separation constraints within the mathematical model.

    \item \textbf{Aircraft Geometry and Buffer:} Each aircraft is represented by its rectangular footprint ($W_a, L_a$). A mandatory safety `Buffer` must be maintained between all aircraft and, crucially, between aircraft and the hangar walls. This dual application of the buffer is vital: the aircraft-to-aircraft buffer ensures safe separation, while the aircraft-to-wall buffer provides essential clearance for personnel and ground-based maintenance equipment (e.g., mobile lifts, hydraulic platforms), ensuring safe maintenance access while mitigating the risk of collisions and scratches, a factor identified as a critical source of danger in recent hangar scheduling models \citep{liuet2023scheduling}. In contrast to more complex geometric approaches like the no-fit polygon (NFP) method used in some literature \citep{qin2018parking, qin2019}, this simplification is justified by operational reality. The space enforced by a rectangular bounding box is practically necessary to accommodate maintenance equipment, allow for tasks such as opening engine cowlings, and ensure safe passage for technicians. This pragmatic choice aligns with other recent models in the field, such as the heuristic framework developed by \citet{zhou2023membership}, which also utilizes buffered rectangles to ensure safe and realistic spatial allocation. While NFP could theoretically achieve a denser packing, the use of buffered rectangles provides a realistic and computationally tractable representation of the required operational space.

    \item \textbf{Static Placement:} Once an aircraft is parked in the hangar, its position is considered fixed until its scheduled departure. The model does not account for the repositioning or shuffling of aircraft during their maintenance stay. This assumption significantly simplifies the model but also aligns with operational practice, where repositioning an aircraft is a complex, time-consuming, and costly procedure that is generally avoided. By creating a stable plan, this assumption enhances predictability. Incorporating dynamic repositioning capabilities remains a viable and interesting avenue for future research.
\end{itemize}

\begin{figure}[htbp!]
    \centering
    
    \begin{subfigure}[b]{0.25\textwidth}
        \centering
        \includegraphics[width=\linewidth]{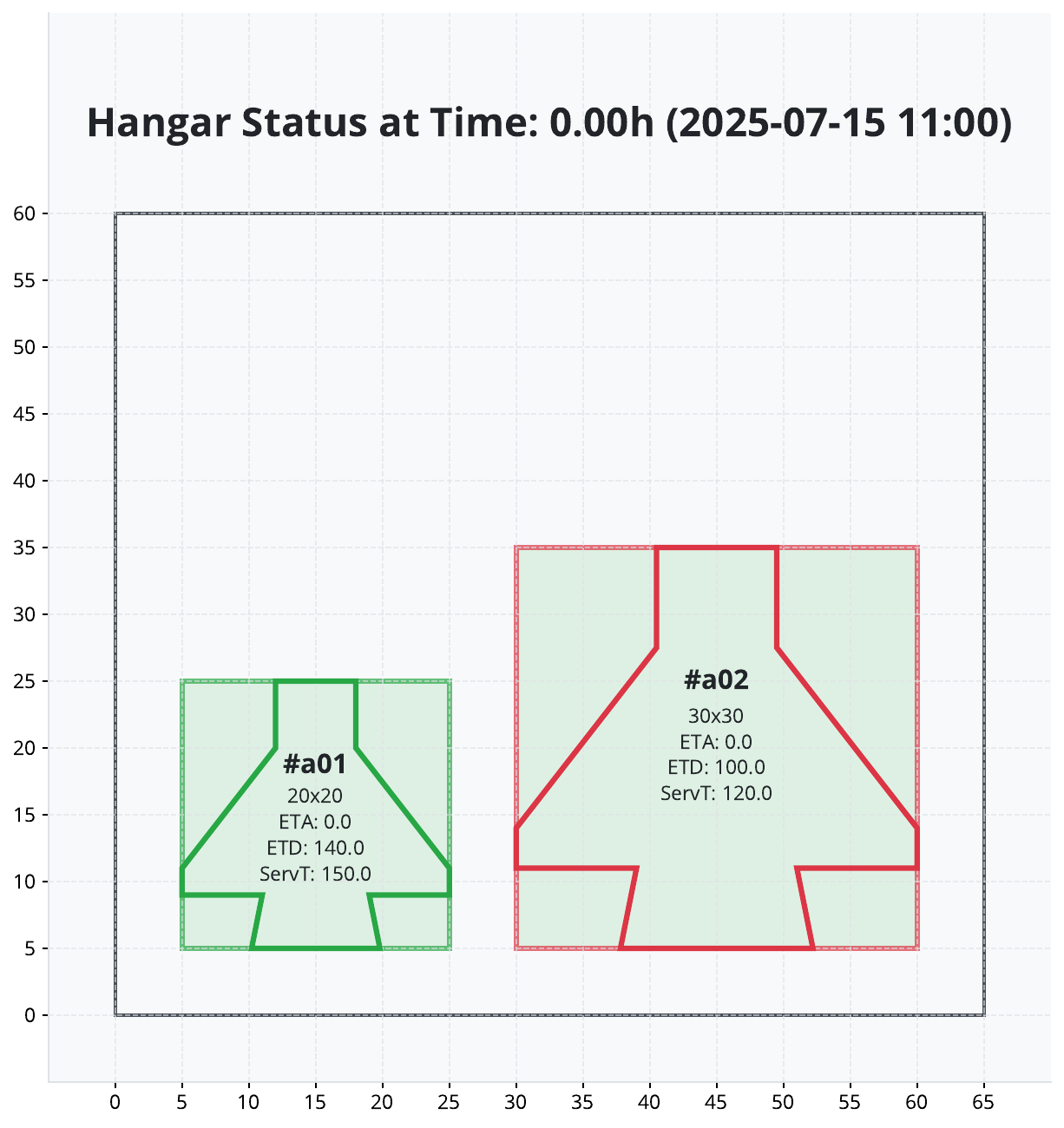}
        \caption{t=0.00h}
        \label{fig:full_state_01}
    \end{subfigure}%
    \begin{subfigure}[b]{0.25\textwidth}
        \centering
        \includegraphics[width=\linewidth]{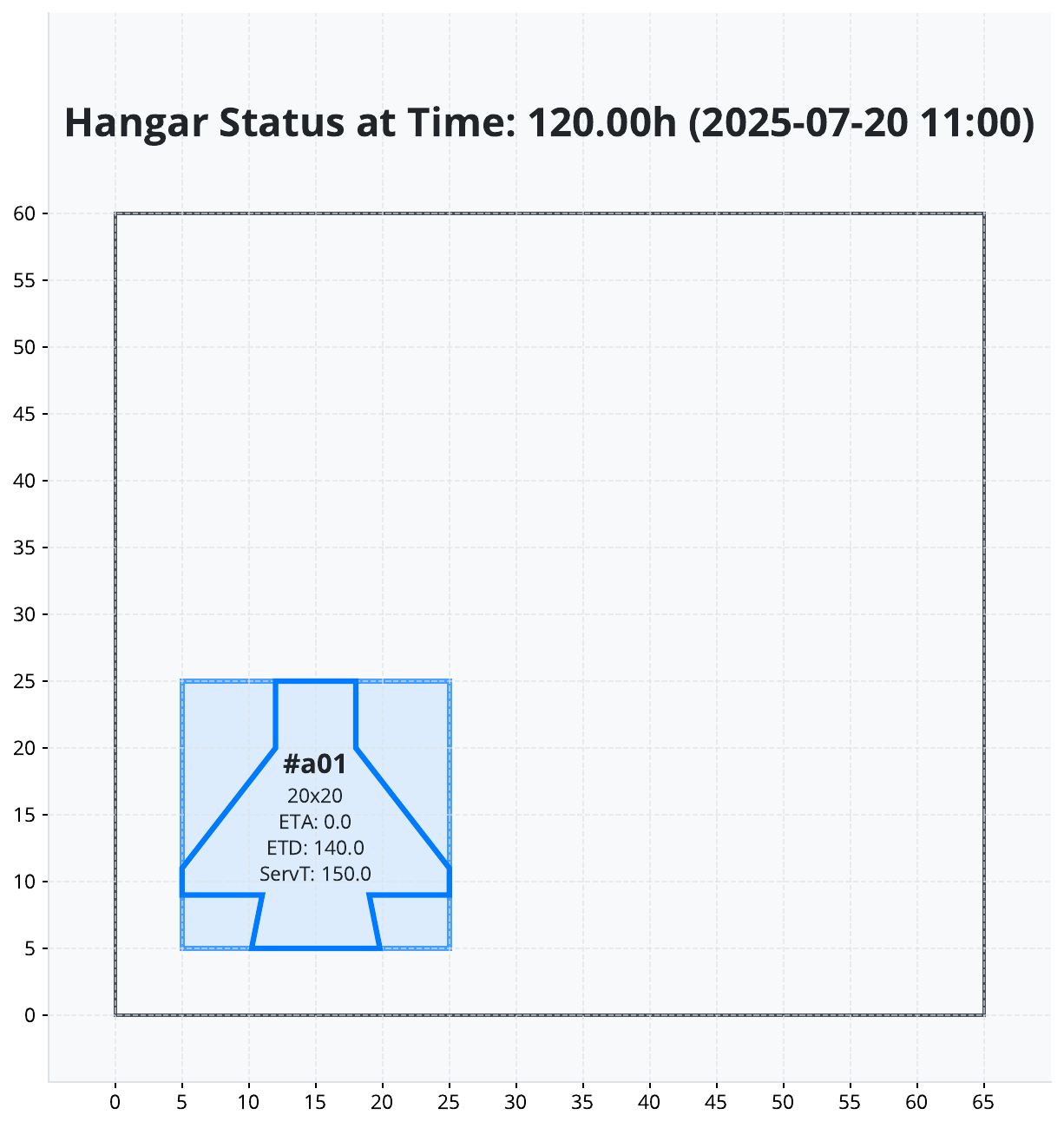}
        \caption{t=120.00h}
        \label{fig:full_state_02}
    \end{subfigure}%
    \begin{subfigure}[b]{0.25\textwidth}
        \centering
        \includegraphics[width=\linewidth]{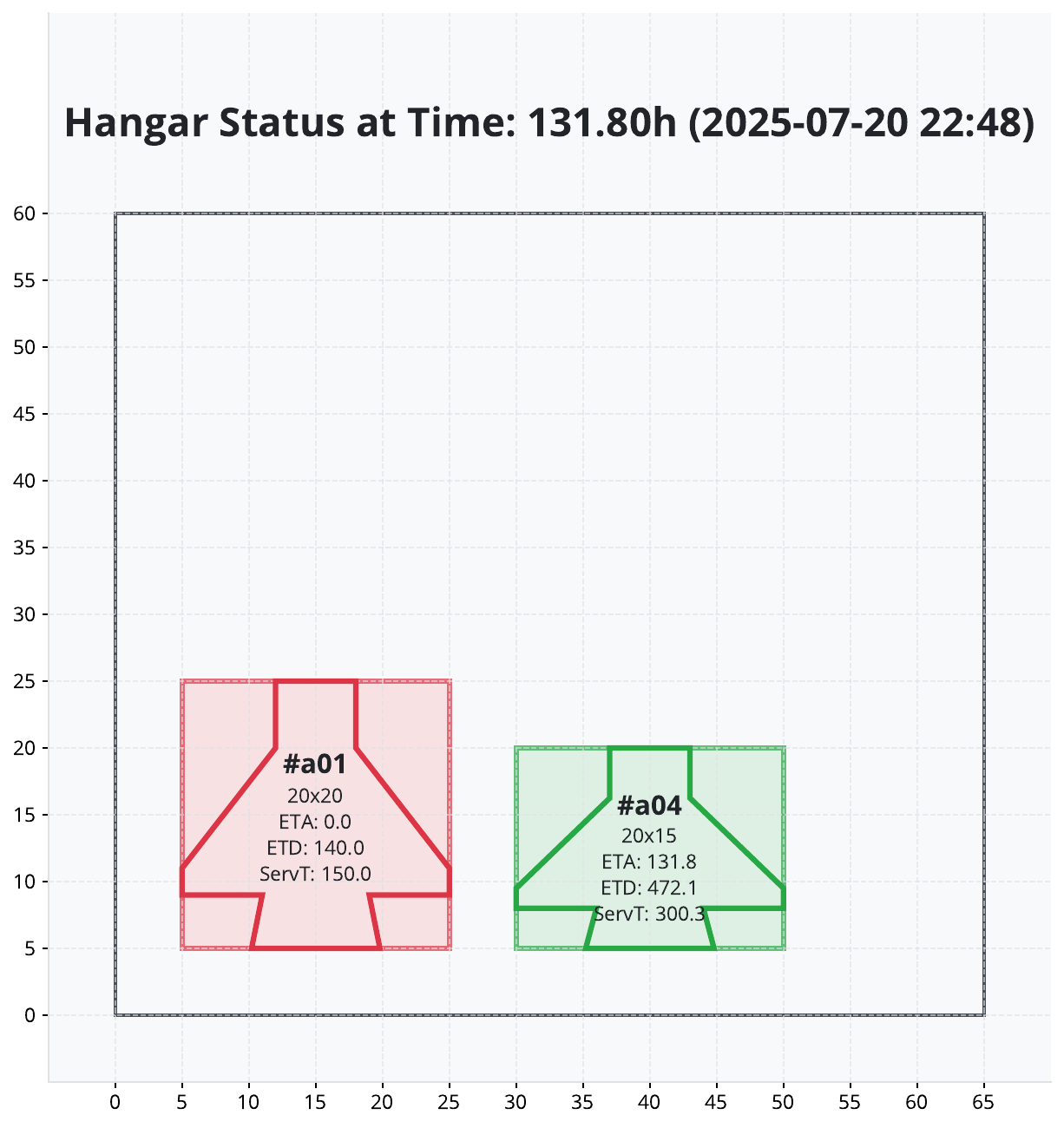}
        \caption{t=131.80h}
        \label{fig:full_state_03}
    \end{subfigure}%
    \begin{subfigure}[b]{0.25\textwidth}
        \centering
        \includegraphics[width=\linewidth]{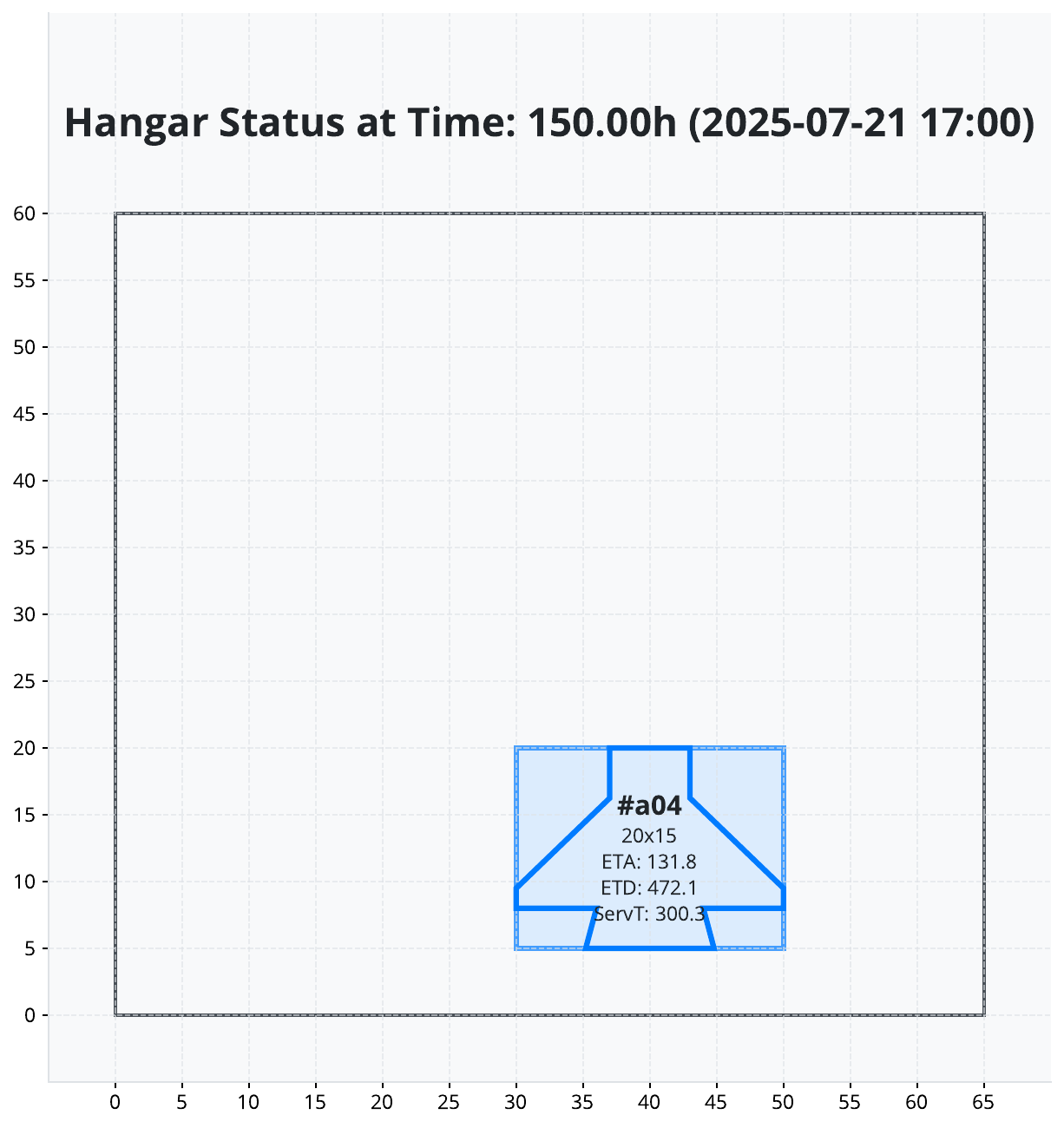}
        \caption{t=150.00h}
        \label{fig:full_state_04}
    \end{subfigure}
    
    \vspace{0.1cm}
    \begin{subfigure}[b]{0.25\textwidth}
        \centering
        \includegraphics[width=\linewidth]{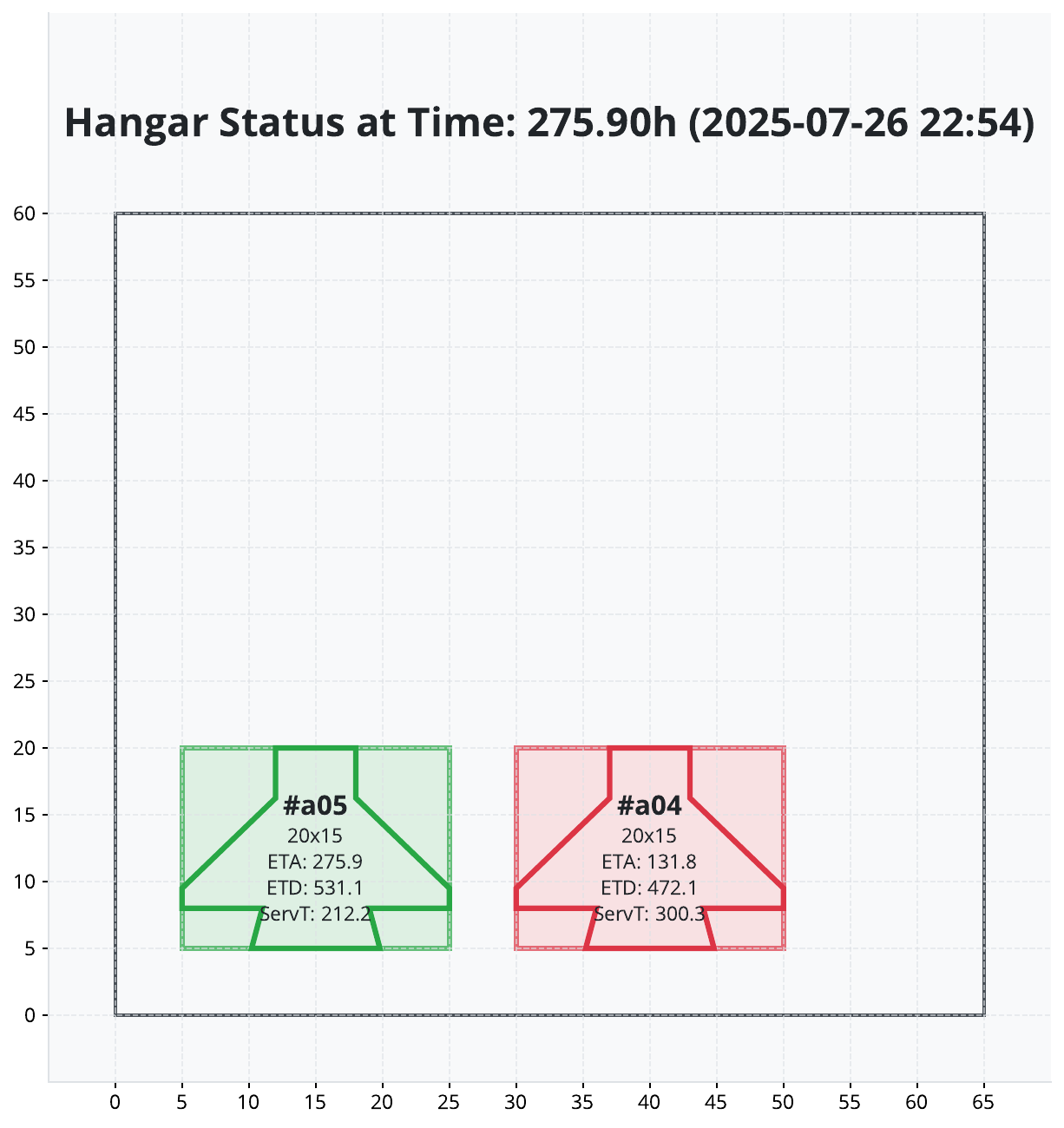}
        \caption{t=275.90h}
        \label{fig:full_state_05}
    \end{subfigure}%
    \begin{subfigure}[b]{0.25\textwidth}
        \centering
        \includegraphics[width=\linewidth]{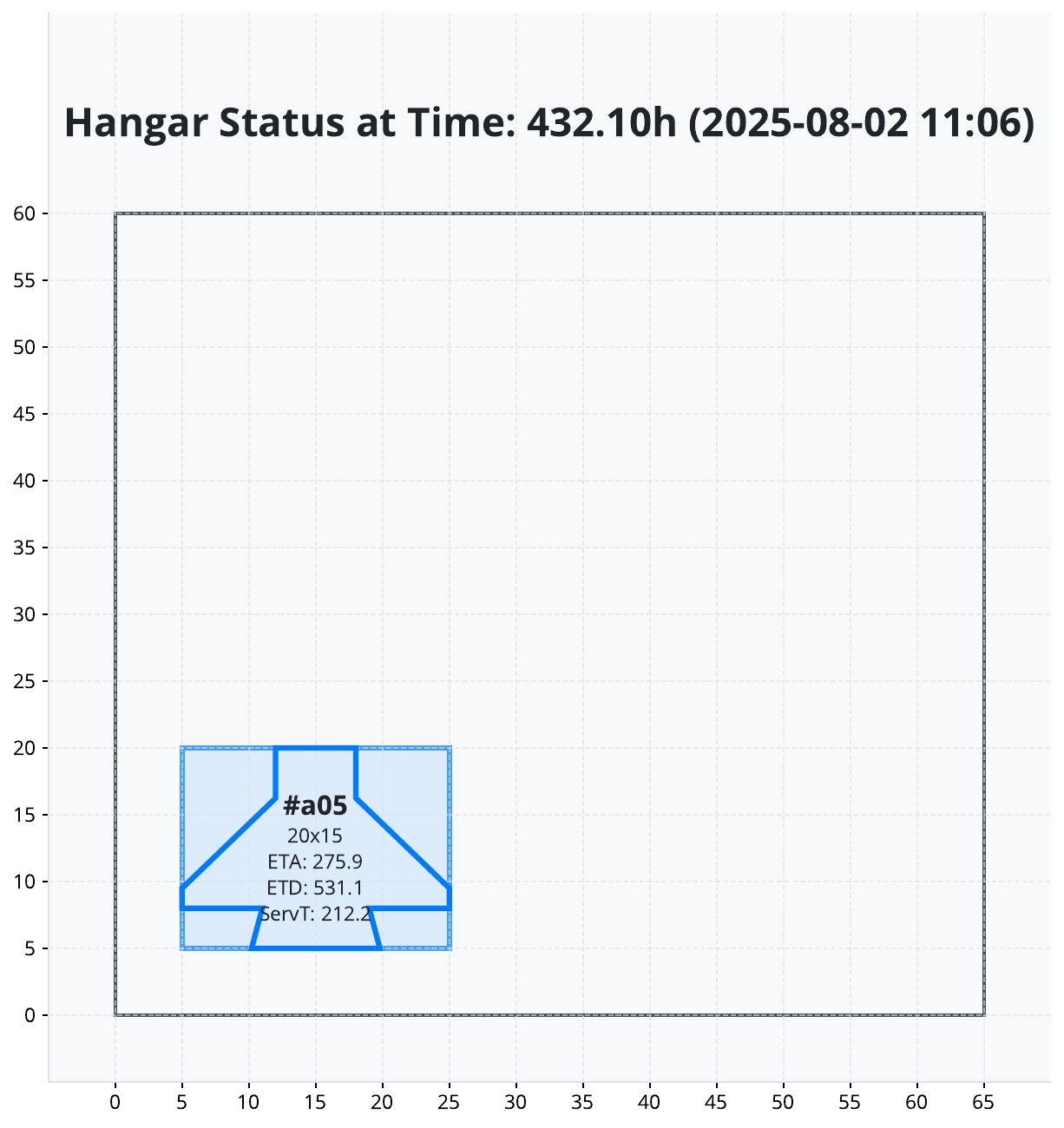}
        \caption{t=432.10h}
        \label{fig:full_state_06}
    \end{subfigure}%
    \begin{subfigure}[b]{0.25\textwidth}
        \centering
        \includegraphics[width=\linewidth]{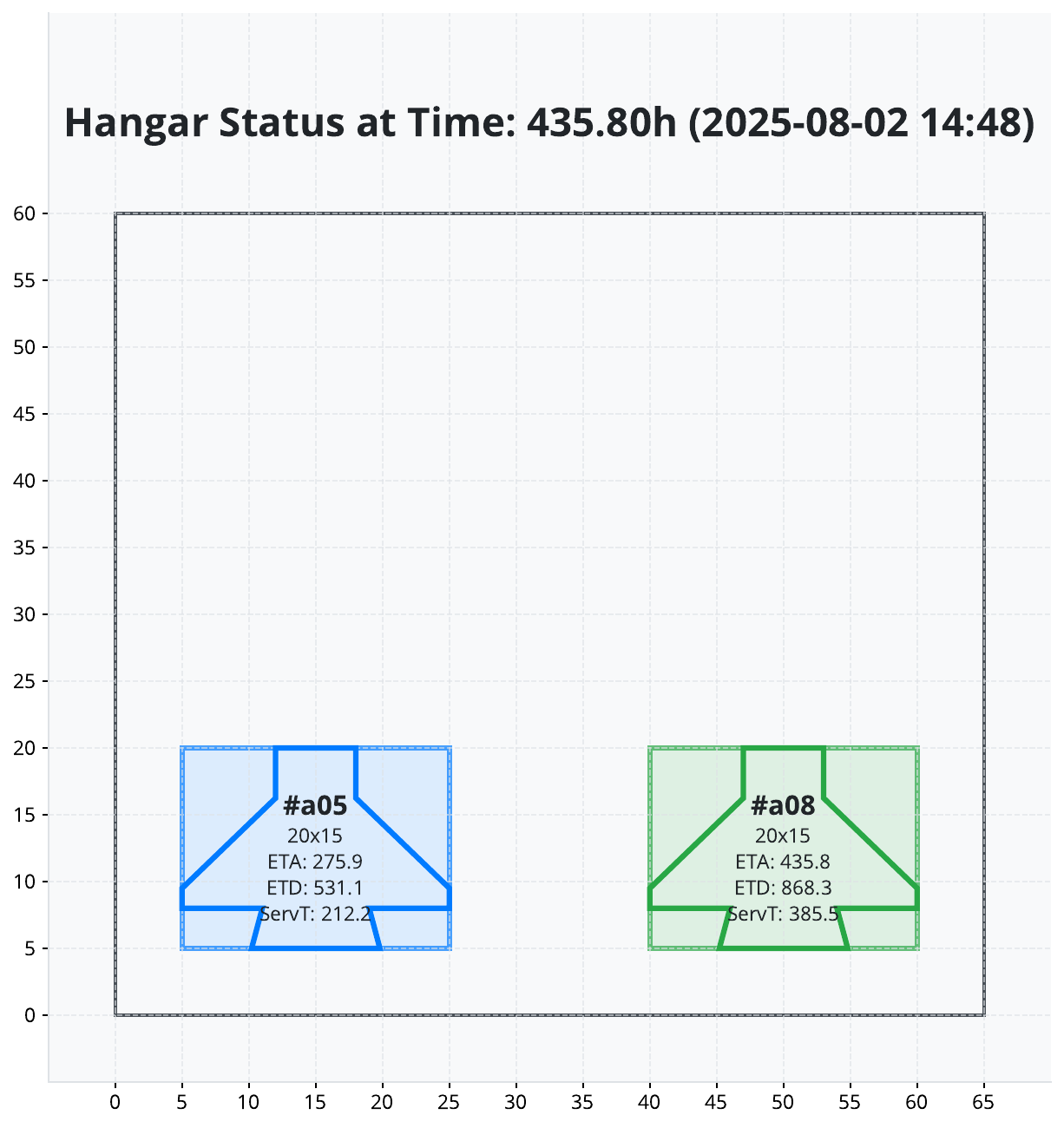}
        \caption{t=435.80h}
        \label{fig:full_state_07}
    \end{subfigure}%
    \begin{subfigure}[b]{0.25\textwidth}
        \centering
        \includegraphics[width=\linewidth]{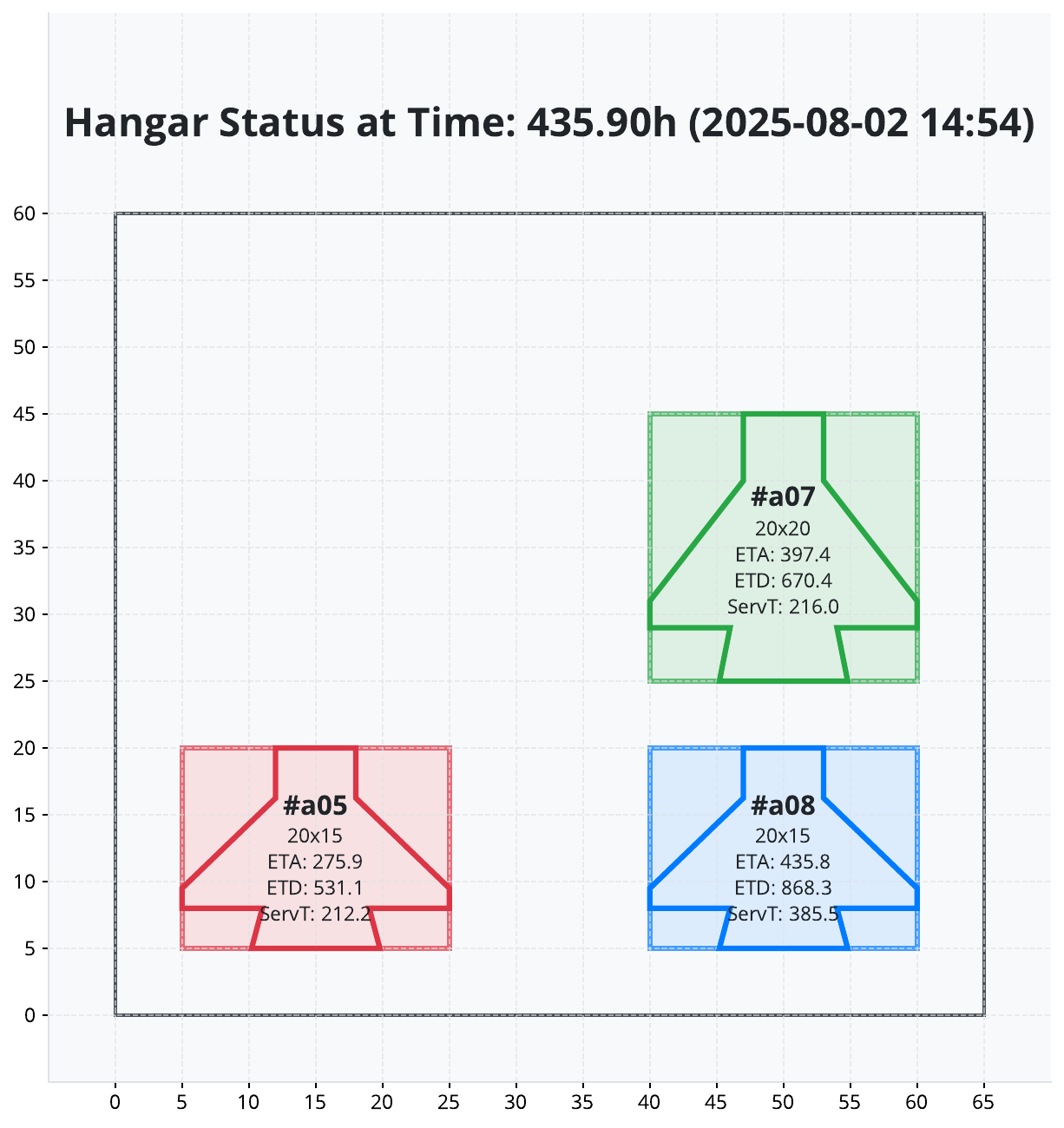}
        \caption{t=435.90h}
        \label{fig:full_state_08}
    \end{subfigure}

    \vspace{0.1cm}
    \begin{subfigure}[b]{0.25\textwidth}
        \centering
        \includegraphics[width=\linewidth]{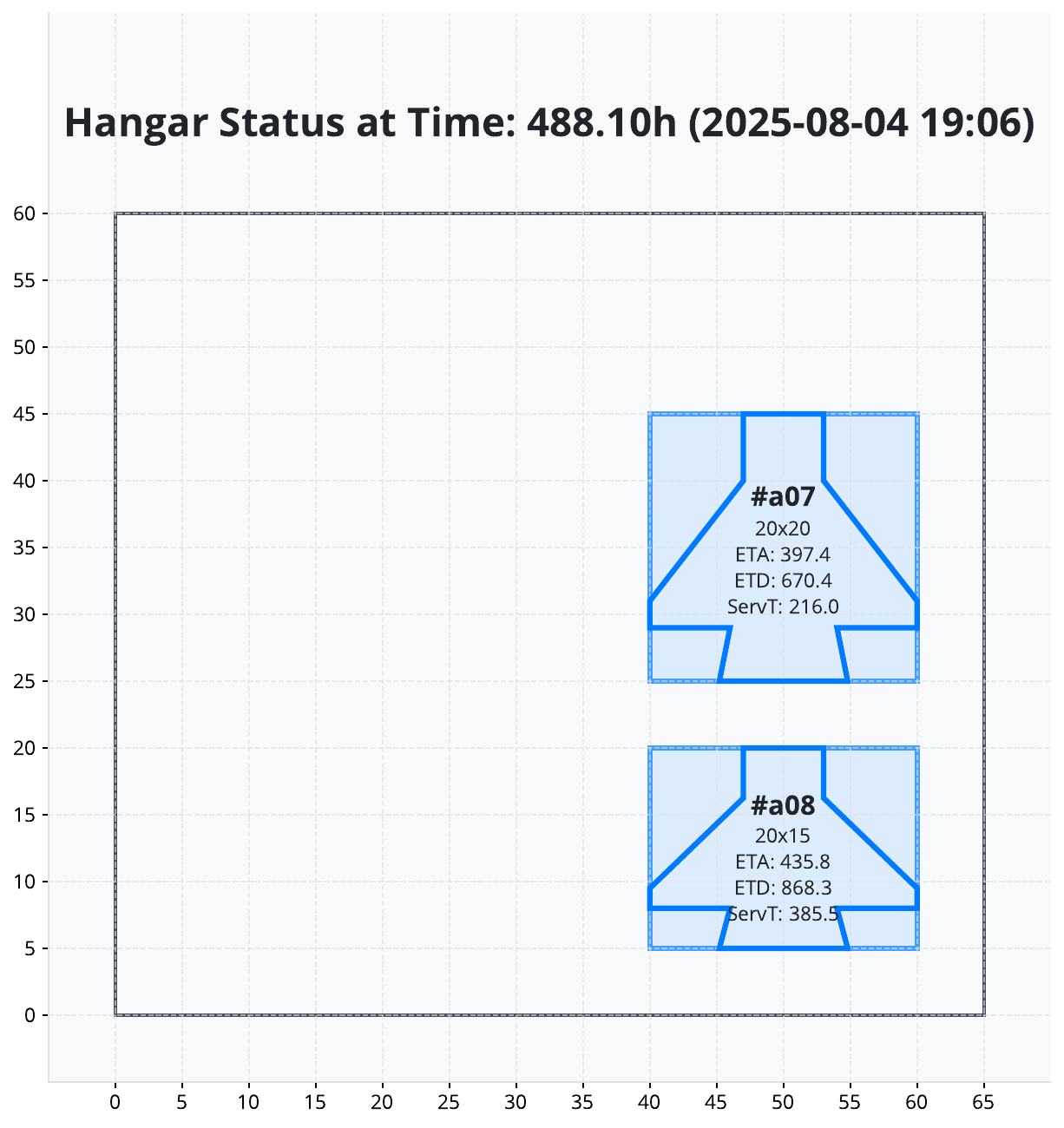}
        \caption{t=488.10h}
        \label{fig:full_state_09}
    \end{subfigure}%
    \begin{subfigure}[b]{0.25\textwidth}
        \centering
        \includegraphics[width=\linewidth]{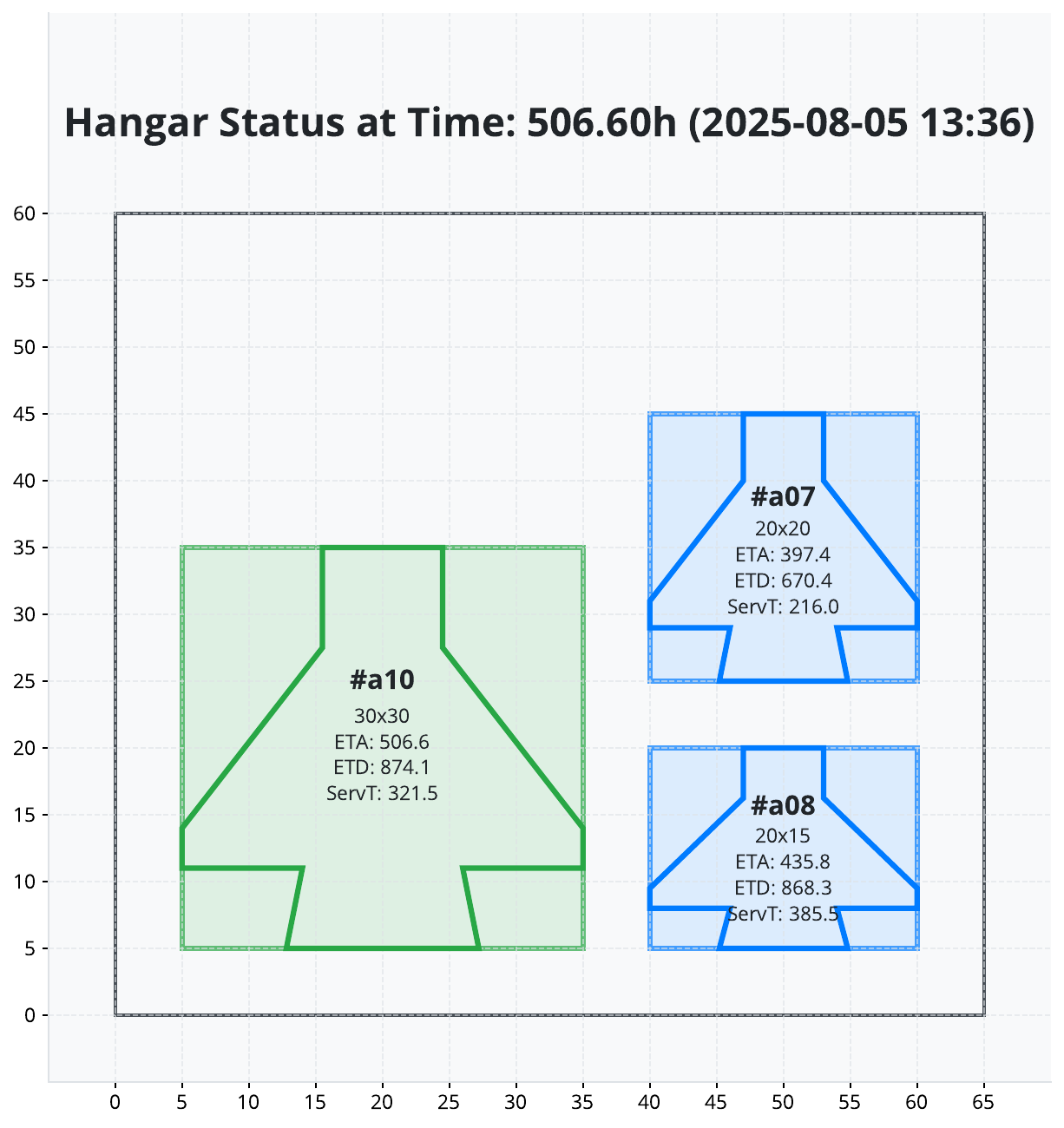}
        \caption{t=506.60h}
        \label{fig:full_state_10}
    \end{subfigure}%
    \begin{subfigure}[b]{0.25\textwidth}
        \centering
        \includegraphics[width=\linewidth]{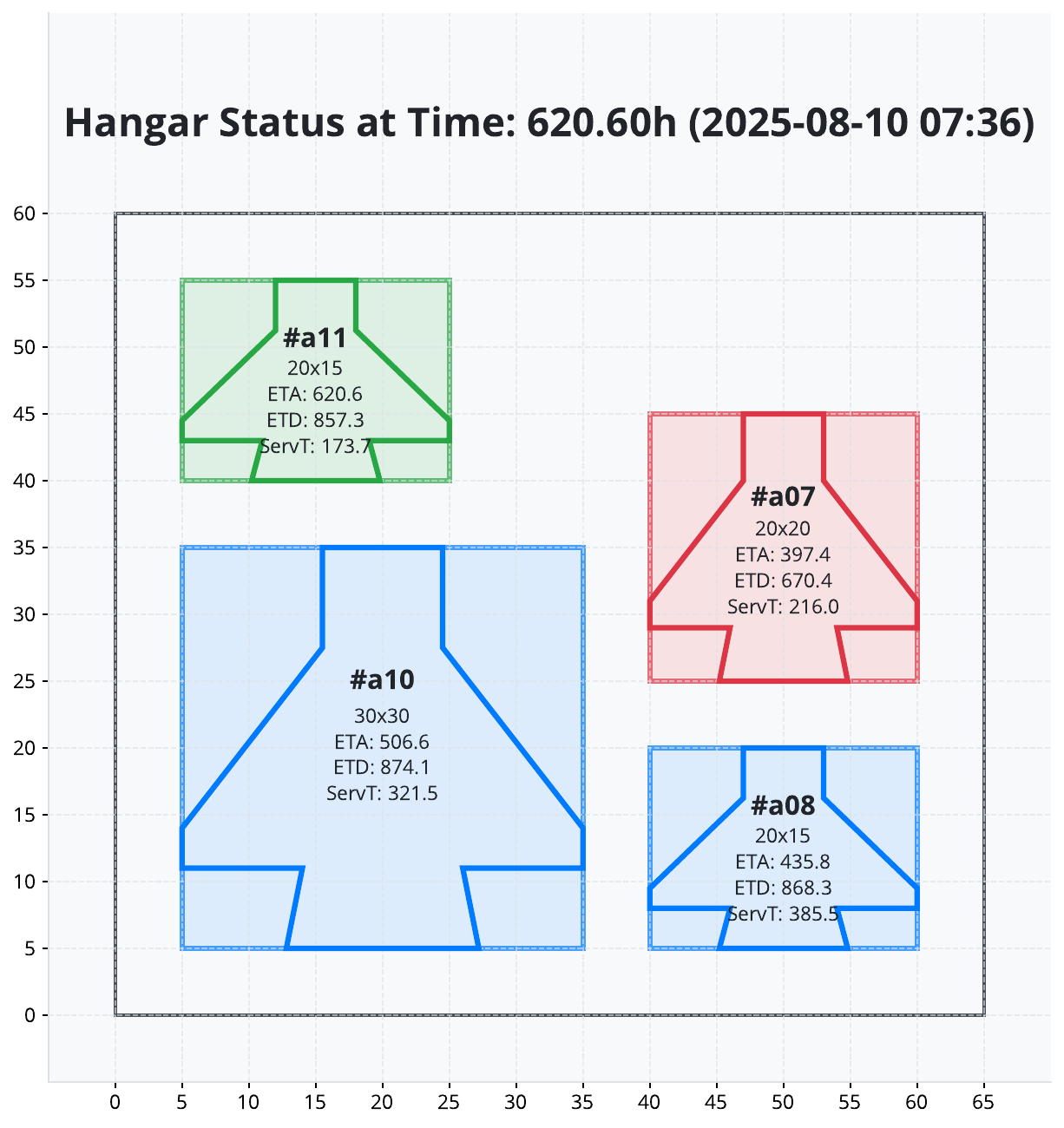}
        \caption{t=620.60h}
        \label{fig:full_state_11}
    \end{subfigure}%
    \begin{subfigure}[b]{0.25\textwidth}
        \centering
        \includegraphics[width=\linewidth]{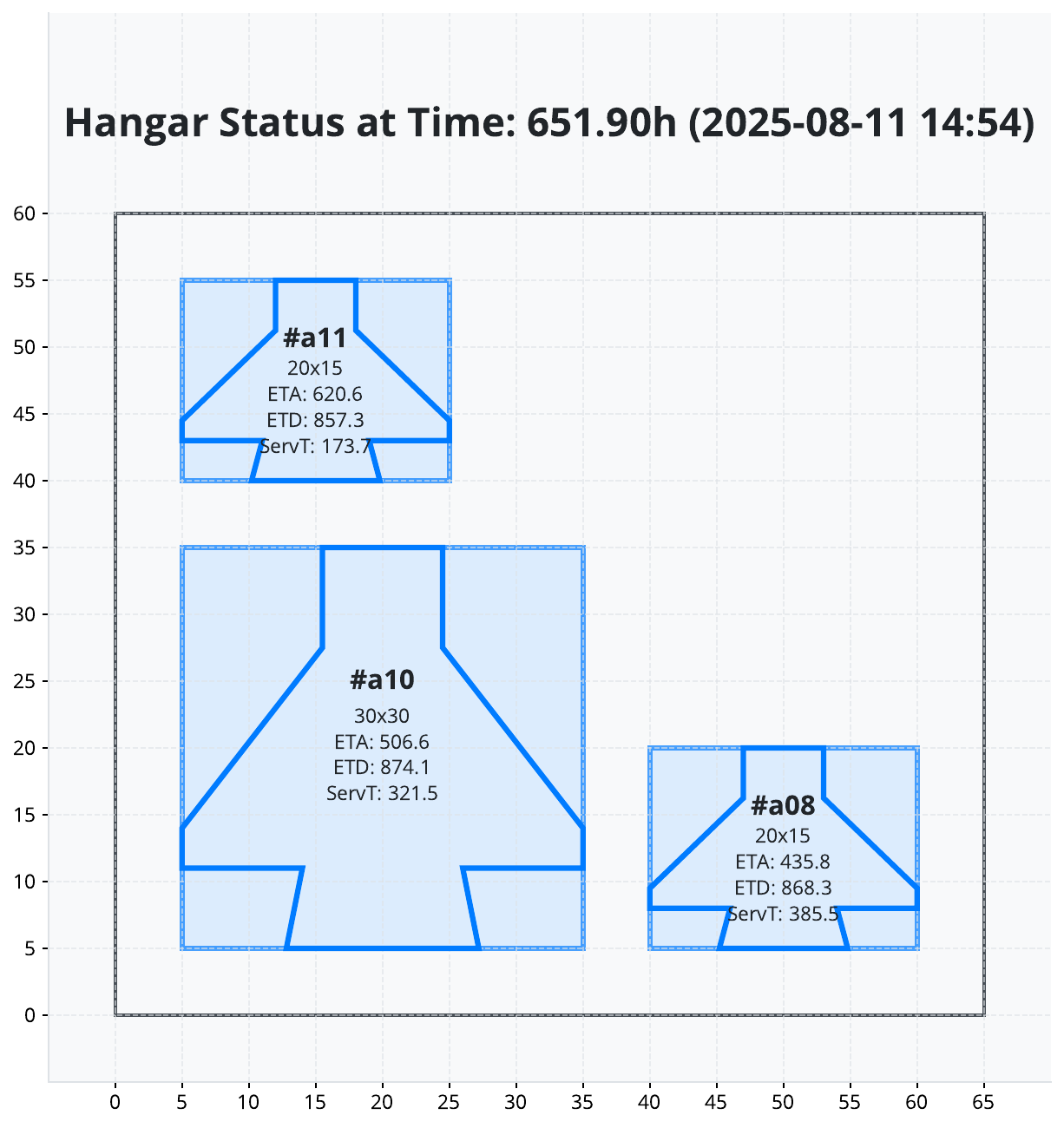}
        \caption{t=651.90h}
        \label{fig:full_state_12}
    \end{subfigure}

    \vspace{0.1cm}
    \begin{subfigure}[b]{0.25\textwidth}
        \centering
        \includegraphics[width=\linewidth]{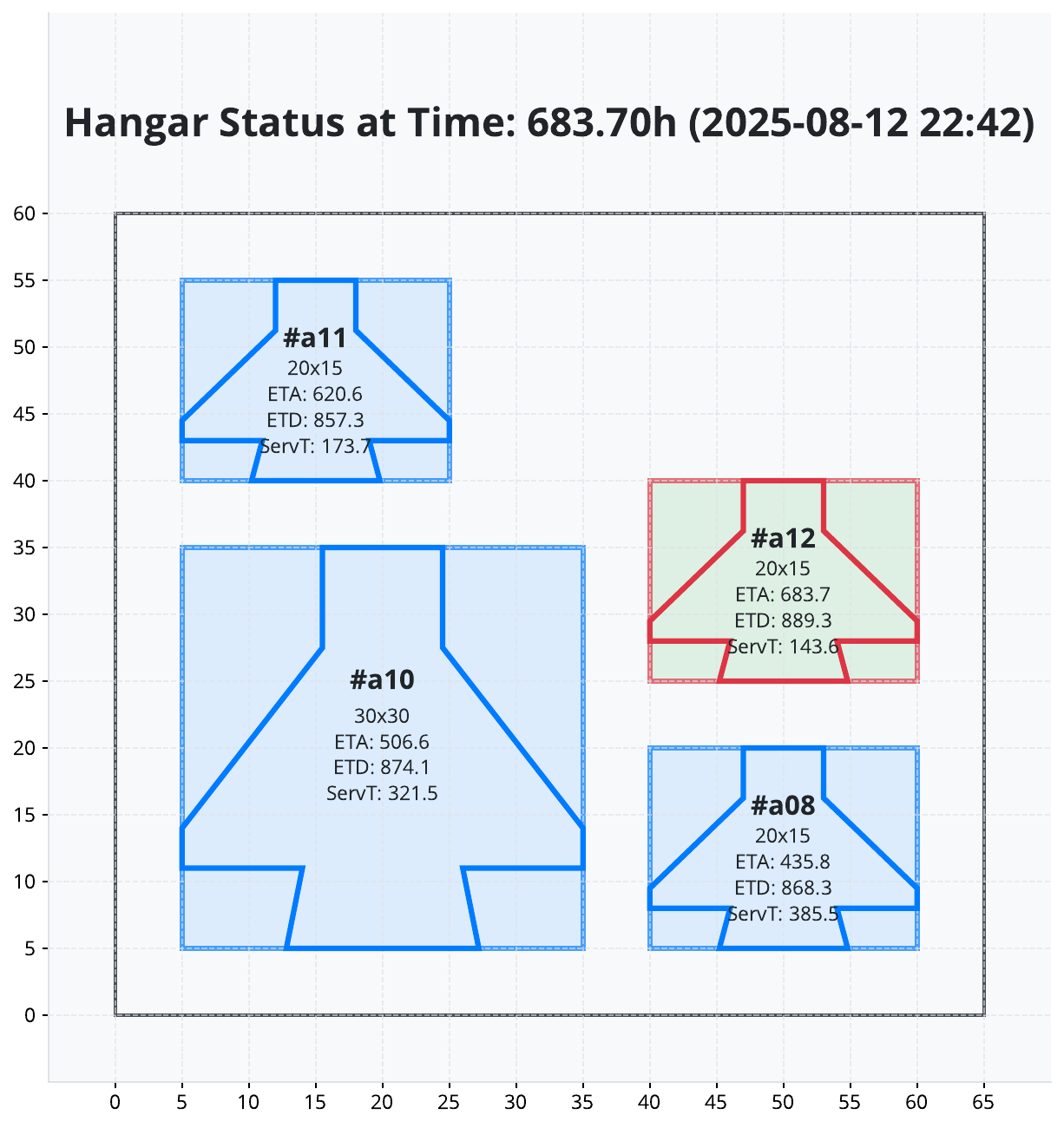}
        \caption{t=683.70h}
        \label{fig:full_state_13}
    \end{subfigure}%
    \begin{subfigure}[b]{0.25\textwidth}
        \centering
        \includegraphics[width=\linewidth]{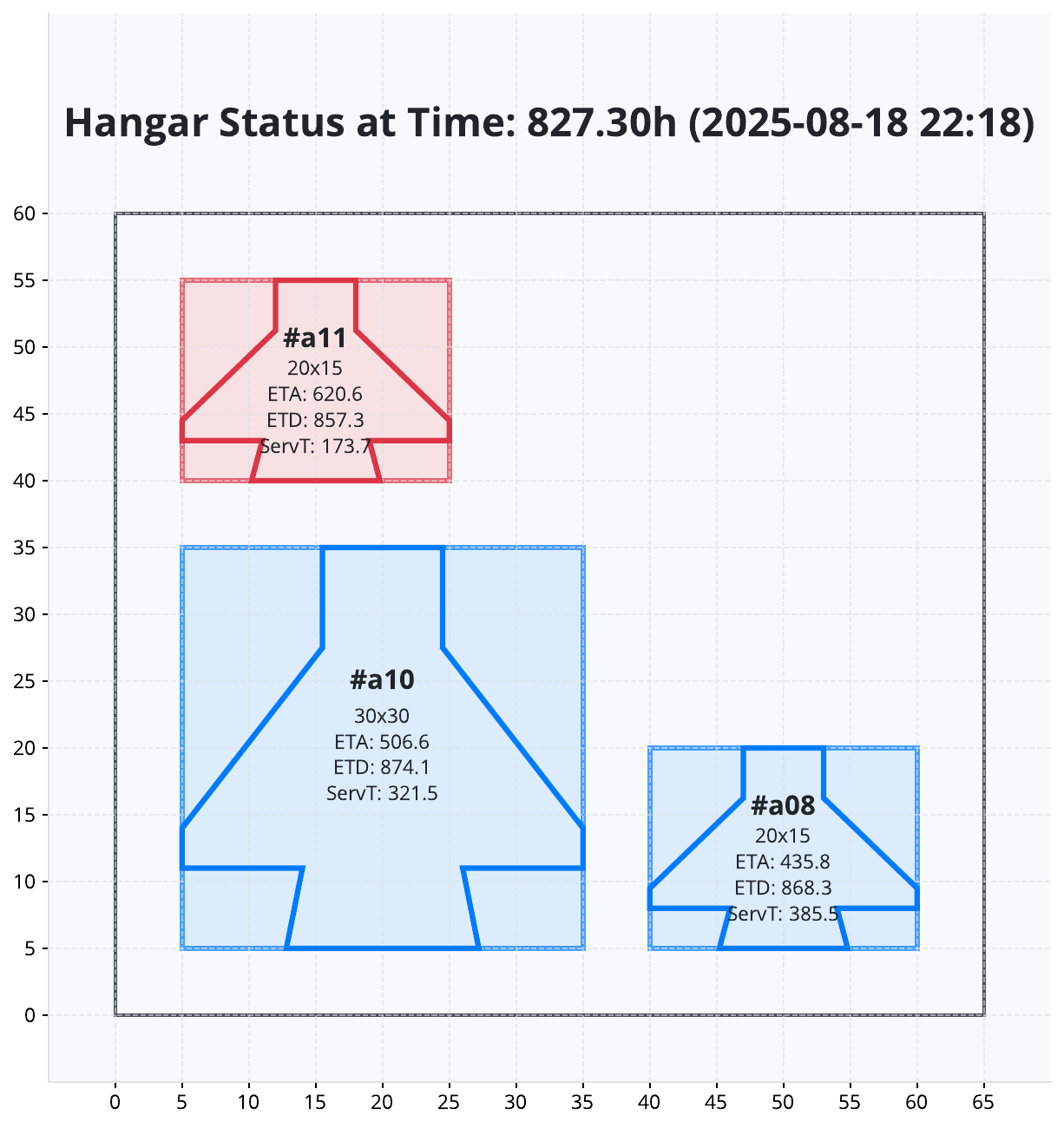}
        \caption{t=827.30h}
        \label{fig:full_state_14}
    \end{subfigure}%
    \begin{subfigure}[b]{0.25\textwidth}
        \centering
        \includegraphics[width=\linewidth]{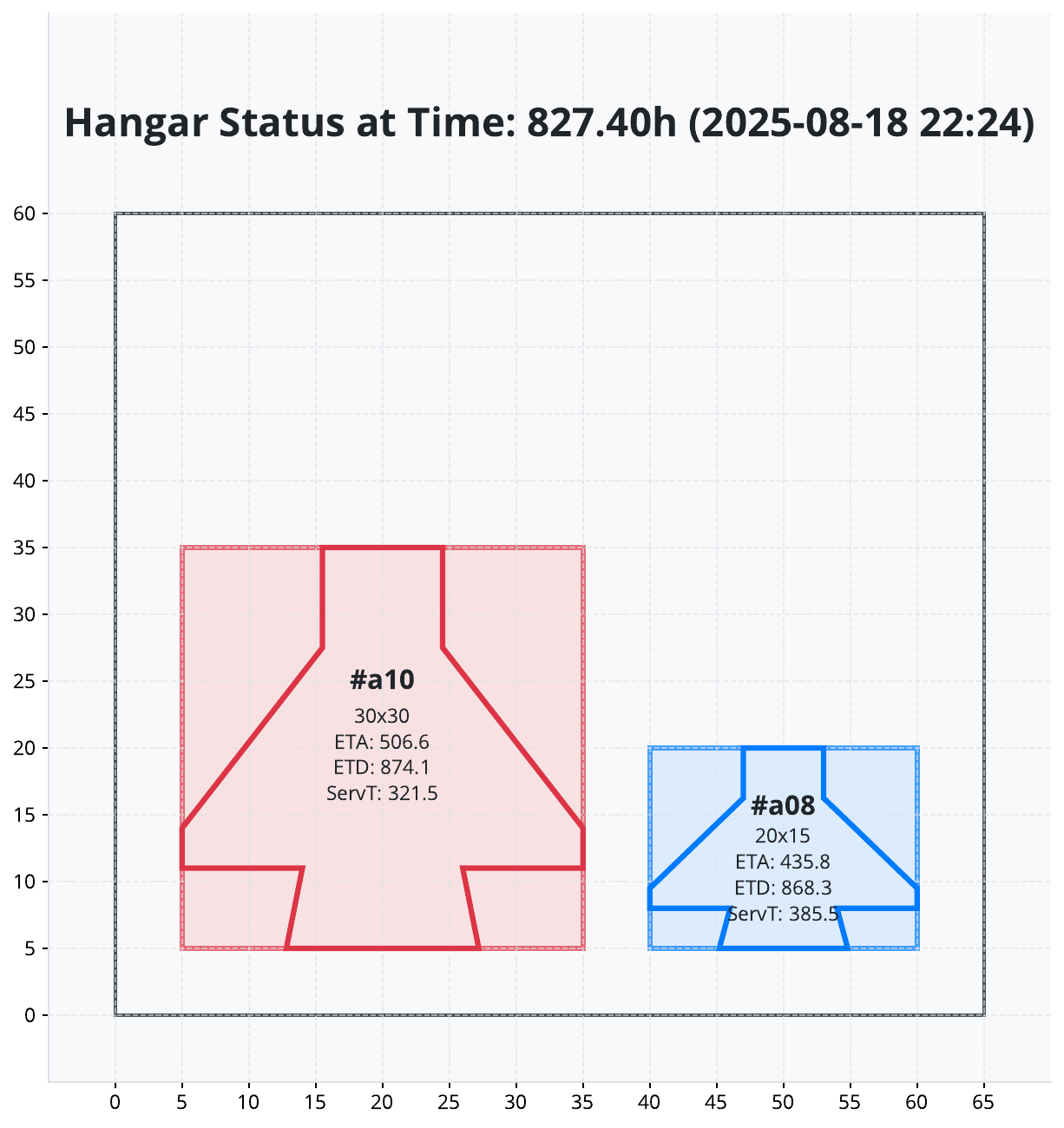}
        \caption{t=827.40h}
        \label{fig:full_state_15}
    \end{subfigure}%
    \begin{subfigure}[b]{0.25\textwidth}
        \centering
        \includegraphics[width=\linewidth]{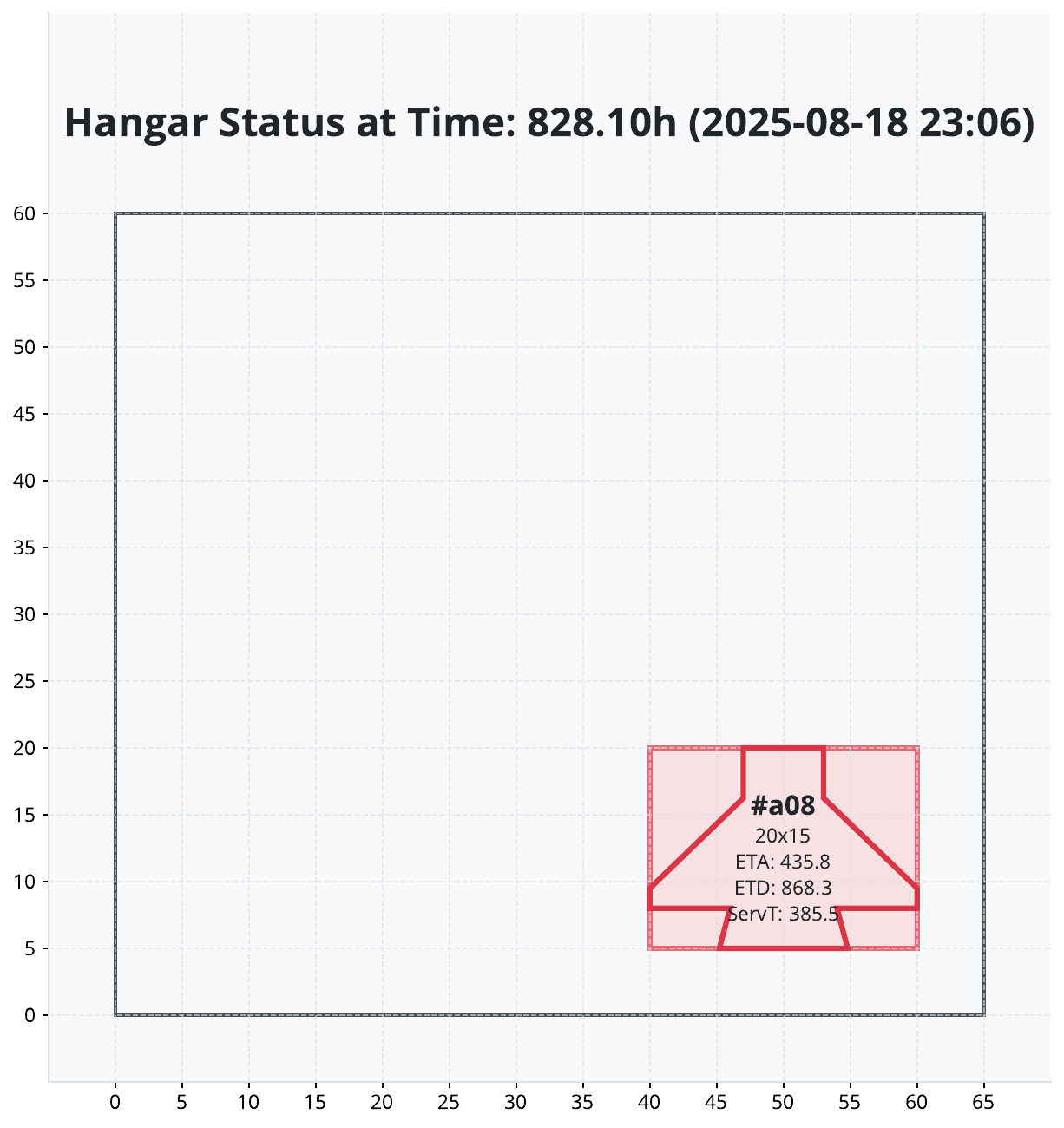}
        \caption{t=828.10h}
        \label{fig:full_state_16}
    \end{subfigure}
    
    \caption{An illustrative optimal solution for a 12-aircraft instance, showcasing the model's ability to manage complex spatial-temporal interactions and blocking constraints. The progression shows the hangar's state at key arrival and departure events. Aircraft colors denote their status: static (blue), arriving (green), or pre-departure (red).}
\label{fig:full_hangar_evolution}
\end{figure}

Figure~\ref{fig:full_hangar_evolution} presents the complete temporal evolution of a hangar configuration for an illustrative instance, showcasing the complex interplay of arrivals, departures, and spatial arrangements that a feasible solution must manage.

\subsection{Sets and Indices}
Sets are used to group the main entities in our model. We define three primary sets to categorize the aircraft.
\begin{itemize}
    \item $a, b \in A$: This is the universal set of all aircraft involved in the planning horizon. It is the union of aircraft already in the hangar and new aircraft requesting service ($A = C \cup F$).
    \item $c, d \in C$: This set represents the \textbf{current} aircraft that are already inside the hangar at the start of the planning period ($C \subset A$). These aircraft have fixed initial positions and remaining service times.
    \item $f, g \in F$: This set includes all \textbf{future} (new) aircraft that have requested hangar space for maintenance ($F \subset A$). The model will decide whether to accept or reject these requests.
\end{itemize}

\subsection{Scalars}
Scalars are constant values that define the physical and operational environment of the model.
\begin{itemize}
    \item $HW, HL$: These represent the total interior \textbf{width} and \textbf{length} of the hangar, respectively. They define the boundaries of the available parking area.
    \item $Buffer$: This scalar specifies the minimum required safety distance that must be maintained between any two aircraft, and between an aircraft and the hangar walls. This is crucial for safe movement and operations.
    \item $M_T, M_X, M_Y$: These are sufficiently large positive constants, commonly known as "Big-M" values, used in mathematical programming to activate or deactivate constraints. Their values are set dynamically based on the problem's dimensions to ensure they are sufficiently large for any given instance without being unnecessarily loose. This practice of generating tight, instance-specific bounds is crucial for preserving the model's computational performance and scalability.
    \begin{itemize}[label=$\circ$]
        \item $M_T = \max_{f \in F}(ETA_f) + \sum_{a \in A}(ServT_a)$: A large time constant, calculated as the latest arrival time plus the sum of all service times. It is guaranteed to be larger than any possible roll-in or roll-out time in the model.
        \item $M_X = HW$: A large distance constant equal to the hangar width.
        \item $M_Y = HL$: A large distance constant equal to the hangar length.
    \end{itemize}
    \item $\varepsilon_t$: A small time value representing the minimum required interval between consecutive aircraft movements (arrivals or departures) to prevent logistical conflicts.
    \item $\varepsilon_p$: A small penalty coefficient used in the objective function to encourage the model to place new aircraft closer to the origin (0,0), promoting a tidy and consolidated layout.
\end{itemize}

\subsection{Parameters}
Parameters are input data specific to each aircraft.
\begin{itemize}
    \item $W_a, L_a$: The \textbf{width} and \textbf{length} of aircraft $a$. These define the rectangular footprint of each aircraft.
    \item $ETA_a, ETD_a$: The \textbf{Expected Time of Arrival} and \textbf{Expected Time of Departure} for aircraft $a$. These are the target times based on the initial schedule.
    \item $ServT_a$: The required \textbf{service time} for aircraft $a$. For current aircraft ($c \in C$), this is the \emph{remaining} service time.
    \item $P^{Rej}_f, P^{Arr}_f, P^{Dep}_a$: These are the penalty costs. $P^{Rej}_f$ is the cost incurred if we \textbf{reject} a new aircraft $f$. $P^{Arr}_f$ is the per-unit-time cost for an \textbf{arrival delay} of aircraft $f$. $P^{Dep}_a$ is the per-unit-time cost for a \textbf{departure delay} of aircraft $a$.
    \item $X^{init}_c, Y^{init}_c$: The fixed initial X and Y coordinates of current aircraft $c$, which are already in the hangar.
\end{itemize}

\subsection{Decision Variables}
Decision variables are the outputs of the model; their values are determined by the solver to achieve the optimal solution.

\subsubsection{Continuous Variables}
\begin{itemize}
    \item $X_a, Y_a$: The X and Y coordinates of the front-left corner of aircraft $a$. These variables determine the exact placement of each aircraft in the hangar.
    \item $Roll\_in_a, Roll\_out_a$: The actual \textbf{roll-in} (entry) and \textbf{roll-out} (exit) times for aircraft $a$.
    \item $D^{Arr}_f, D^{Dep}_a$: The calculated \textbf{arrival delay} for new aircraft $f$ and \textbf{departure delay} for any aircraft $a$.
\end{itemize}
\subsubsection{Binary Variables}
\begin{itemize}
    \item $Accept_a$: A binary flag that is 1 if aircraft $a$ is accepted for service, and 0 otherwise. For current aircraft ($c \in C$), this is fixed to 1.
    \item $Right_{ab}$: 1 if aircraft $a$ is positioned completely to the \textbf{right} of aircraft $b$, 0 otherwise.
    \item $Above_{ab}$: 1 if aircraft $a$ is positioned completely \textbf{above} (in the direction of increasing Y) aircraft $b$, 0 otherwise.
    \item $OutIn_{ab}$: 1 if aircraft $a$ rolls \textbf{out} before aircraft $b$ rolls \textbf{in}.
    \item $InIn_{ab}$: 1 if aircraft $a$ rolls \textbf{in} before aircraft $b$.
    \item $OutOut_{ab}$: 1 if aircraft $a$ rolls \textbf{out} before aircraft $b$.
    \item $InOut_{ab}$: 1 if aircraft $a$ rolls \textbf{in} before aircraft $b$ rolls out.
\end{itemize}

\subsection{Objective Function}
The objective is to minimize the total operational cost $Z$, which is formulated as a weighted sum of four distinct cost components.
\begin{equation}
\begin{aligned}
\min Z = & \sum_{f \in F} P^{Rej}_f(1 - Accept_f) + \sum_{f \in F} P^{Arr}_f D^{Arr}_f \\
& + \sum_{a \in A} P^{Dep}_a D^{Dep}_a + \varepsilon_p \sum_{f \in F} (X_f + Y_f)
\end{aligned}
\end{equation}
\textbf{Explanation:}
\begin{itemize}
    \item \textbf{Rejection Cost:} $\sum_{f \in F} P^{Rej}_f(1 - Accept_f)$. This term imposes a penalty for each new service request $f$ that is rejected. If $Accept_f=0$, the term $(1 - Accept_f)$ becomes 1, and the cost $P^{Rej}_f$ is incurred.
    \item \textbf{Arrival Delay Cost:} $\sum_{f \in F} P^{Arr}_f D^{Arr}_f$. This component penalizes deviations from the planned schedule by calculating the cost of arrival delays for new aircraft.
    \item \textbf{Departure Delay Cost:} $\sum_{a \in A} P^{Dep}_a D^{Dep}_a$. Similarly, this term accounts for the cost associated with departure delays for all aircraft, both current and new.
    \item \textbf{Positioning Cost:} $\varepsilon_p \sum_{f \in F} (X_f + Y_f)$. This is a regularization term with a small coefficient $\varepsilon_p$ that encourages the model to place new aircraft closer to the hangar's origin. This promotes a compact and organized spatial layout without overriding the primary cost considerations.
\end{itemize}

\subsubsection{Objective Function Rationale and Equivalence to Profit Maximization}
For modeling convenience, our objective function is formulated to minimize total costs. However, it is crucial to interpret the "rejection penalty," \texttt{$P^{Rej}_f$}, not as a literal penalty, but as the \textbf{opportunity cost} or \textbf{expected profit} lost by rejecting a service request for aircraft \textit{f}. This aligns the model's logic with the real-world goal of maximizing profitability.

To formally demonstrate this, we can show that our cost-minimization objective is mathematically equivalent to a profit-maximization objective. Let the total delay cost be defined as $D = \sum_{f \in F} P^{Arr}_f D^{Arr}_f + \sum_{a \in A} P^{Dep}_a D^{Dep}_a$.

Our objective is to minimize total cost, $Z_{\text{cost}}$:
\begin{equation*}
    \min Z_{\text{cost}} = \sum_{f \in F} P^{Rej}_f (1 - Accept_f) + D
\end{equation*}

An equivalent objective would be to maximize total profit, $Z_{\text{profit}}$, where we gain profit from accepted aircraft and subtract delay costs:
\begin{equation*}
    \max Z_{\text{profit}} = \sum_{f \in F} P^{Rej}_f \cdot Accept_f - D
\end{equation*}

The equivalence can be proven by rearranging the cost function. Let $C = \sum_{f \in F} P^{Rej}_f$ be the total potential profit if all new requests were accepted. Since $P^{Rej}_f$ are fixed parameters, $C$ is a constant for any given instance.
\begin{align*}
    Z_{\text{cost}} &= \sum_{f \in F} P^{Rej}_f (1 - Accept_f) + D \\
    &= \underbrace{\sum_{f \in F} P^{Rej}_f}_{\text{Constant, } C} - \sum_{f \in F} P^{Rej}_f \cdot Accept_f + D \\
    &= C - \left( \underbrace{\sum_{f \in F} P^{Rej}_f \cdot Accept_f - D}_{Z_{\text{profit}}} \right) \\
    & \implies \boxed{Z_{\text{cost}} = C - Z_{\text{profit}}}
\end{align*}
Since $C$ is a constant, minimizing $Z_{\text{cost}}$ is perfectly equivalent to maximizing $Z_{\text{profit}}$. Therefore, the optimal plan generated by our model is guaranteed to be the most profitable one.

\subsection{Model Constraints}
Constraints define the rules and physical limitations of the system, ensuring that the solution is both feasible and realistic.

\subsubsection{Acceptance and Scheduling Constraints}
These constraints establish the fundamental relationships between the acceptance decision and the scheduling variables.
\begin{align}
& \begin{aligned}
    & X_f + Y_f + Roll\_in_f + Roll\_out_f + D^{Arr}_f + D^{Dep}_f \\
    & \qquad \le (M_X + M_Y + 4M_T) \cdot Accept_f,
  \end{aligned} && \forall f \in F \label{eq:accept} \\
& Roll\_in_f \ge ETA_f \cdot Accept_f, && \forall f \in F \label{eq:rollin_eta} \\
& Roll\_out_a - Roll\_in_a \ge ServT_a \cdot Accept_a, && \forall a \in A \label{eq:servt} \\
& D^{Arr}_f \ge Roll\_in_f - ETA_f, && \forall f \in F \label{eq:darr} \\
& D^{Dep}_a \ge Roll\_out_a - ETD_a, && \forall a \in A \label{eq:ddep}
\end{align}
\textbf{Explanation:}
\begin{itemize}
    \item \textbf{Constraint \eqref{eq:accept}:} A Big-M formulation that ensures if a new aircraft $f$ is rejected ($Accept_f = 0$), all associated continuous variables (position, time, and delay) are forced to be zero. If accepted ($Accept_f = 1$), the constraint becomes non-binding.
    \item \textbf{Constraint \eqref{eq:rollin_eta}:} An accepted aircraft cannot roll in before its Expected Time of Arrival (ETA).
    \item \textbf{Constraint \eqref{eq:servt}:} The duration an accepted aircraft occupies a hangar spot must be sufficient to cover its required service time.
    \item \textbf{Constraint \eqref{eq:darr}:} Defines the arrival delay for a new aircraft as the non-negative difference between its actual roll-in time and its ETA.
    \item \textbf{Constraint \eqref{eq:ddep}:} Defines the departure delay for any aircraft as the non-negative difference between its actual roll-out time and its ETD.
\end{itemize}

\subsubsection{Hangar Physical and Non-overlapping Constraints}
These constraints ensure that aircraft fit within the hangar boundaries and do not physically overlap.
\begin{align}
& X_f \ge Buffer \cdot Accept_f, && \forall f \in F \label{eq:xmin} \\
& X_f + W_f \le HW - Buffer + M_X(1 - Accept_f), && \forall f \in F \label{eq:xmax} \\
& Y_f \ge Buffer \cdot Accept_f, && \forall f \in F \label{eq:ymin} \\
& Y_f + L_f \le HL - Buffer + M_Y(1 - Accept_f), && \forall f \in F \label{eq:ymax} \\
& X_b + W_b + Buffer \le X_a + M_X(1 - Right_{ab}), && \forall a, b \in A, a \neq b \label{eq:right} \\
& Y_b + L_b + Buffer \le Y_a + M_Y(1 - Above_{ab}), && \forall a, b \in A, a \neq b \label{eq:above}
\end{align}
\textbf{Explanation:}
\begin{itemize}
    \item \textbf{Constraints \eqref{eq:xmin} - \eqref{eq:ymax}:} These four constraints enforce that any accepted new aircraft is placed entirely within the hangar's usable area, maintaining the specified safety `Buffer` from the walls.
    \item \textbf{Constraint \eqref{eq:right}:} This constraint, along with its counterpart for the Y-axis, prevents spatial overlap. If $Right_{ab}=1$, it signifies that aircraft $a$ is entirely to the right of aircraft $b$. The constraint then enforces that the X-coordinate of $a$'s left edge ($X_a$) must be greater than or equal to the X-coordinate of $b$'s right edge ($X_b + W_b$) plus the safety buffer. If $Right_{ab}=0$, the Big-M term renders the constraint inactive.
    \item \textbf{Constraint \eqref{eq:above}:} This functions identically to the previous constraint but for the vertical (Y) axis, using the binary variable $Above_{ab}$.
\end{itemize}

\subsubsection{Aircraft Separation and Relationship Constraints}
These constraints establish the logical disjunctive relationships required to avoid conflicts between any pair of aircraft.
\begin{align}
& \begin{aligned}
& Right_{ba} + Right_{ab} + Above_{ba} + Above_{ab} \\
& + OutIn_{ab} + OutIn_{ba} \ge Accept_a + Accept_b - 1,
\end{aligned} && \forall a, b \in A, a < b \label{eq:rel} \\
& Roll\_out_a + \varepsilon_t \le Roll\_in_b + M_T(1 - OutIn_{ab}), && \forall a, b \in A, a \neq b \label{eq:outin} \\
& \begin{aligned}
& Roll\_in_g \ge Roll\_in_f + \varepsilon_t - M_T(1 - InIn_{fg}) \\
& - M_T(2 - Accept_f - Accept_g),
\end{aligned} && \forall f, g \in F, f \neq g \label{eq:inin1} \\
& \begin{aligned}
& Roll\_in_f \ge Roll\_in_g + \varepsilon_t - M_T \cdot InIn_{fg} \\
& - M_T(2 - Accept_f - Accept_g),
\end{aligned} && \forall f, g \in F, f \neq g \label{eq:inin2} \\
& \begin{aligned}
& Roll\_out_b \ge Roll\_out_a + \varepsilon_t - M_T(1 - OutOut_{ab}) \\
& - M_T(2 - Accept_a - Accept_b),
\end{aligned} && \forall a, b \in A, a < b \label{eq:outout1} \\
& \begin{aligned}
& Roll\_out_a \ge Roll\_out_b + \varepsilon_t - M_T \cdot OutOut_{ab} \\
& - M_T(2 - Accept_a - Accept_b),
\end{aligned} && \forall a, b \in A, a < b \label{eq:outout2} \\
& \begin{aligned}
& Roll\_out_b \ge Roll\_in_a + \varepsilon_t - M_T(1 - InOut_{ab}) \\
& - M_T(2 - Accept_a - Accept_b),
\end{aligned} \label{eq:inout1} \\
& \qquad \qquad \qquad \qquad \qquad \forall a, b \in A, a \neq b, \text{ where } a \in F \lor b \in F \nonumber \\
& \begin{aligned}
& Roll\_in_a \ge Roll\_out_b + \varepsilon_t - M_T \cdot InOut_{ab} \\
& - M_T(2 - Accept_a - Accept_b),
\end{aligned} \label{eq:inout2} \\
& \qquad \qquad \qquad \qquad \qquad \forall a, b \in A, a \neq b, \text{ where } a \in F \lor b \in F \nonumber
\end{align}
\textbf{Explanation:}
\begin{itemize}
    \item \textbf{Constraint \eqref{eq:rel}:} This is a generalized disjunctive constraint. It mandates that for any pair of accepted aircraft ($a, b$), they must be separated in at least one dimension: spatially (one is right of, left of, above, or below the other) or temporally (one departs before the other arrives).
    \item \textbf{Constraint \eqref{eq:outin}:} This constraint gives meaning to the $OutIn_{ab}$ variable. If $OutIn_{ab}=1$, it enforces that aircraft $a$ must complete its roll-out at least $\varepsilon_t$ time units before aircraft $b$ can begin its roll-in.
    \item \textbf{Constraints \eqref{eq:inin1} and \eqref{eq:inin2}:} This pair of constraints establishes a strict sequence for the roll-in times of any two new aircraft, $f$ and $g$, if both are accepted. The binary variable $InIn_{fg}$ acts as a switch: if $InIn_{fg}=1$, constraint \eqref{eq:inin2} forces $f$ to roll-in before $g$; if $InIn_{fg}=0$, constraint \eqref{eq:inin1} forces $g$ to roll-in before $f$.
    \item \textbf{Constraints \eqref{eq:outout1} and \eqref{eq:outout2}:} This pair works identically to the one above but enforces a strict sequence for the roll-out times of any two aircraft, preventing simultaneous departures.
    \item \textbf{Constraints \eqref{eq:inout1} and \eqref{eq:inout2}:} This pair establishes a sequential relationship between the roll-in of one aircraft and the roll-out of another, using the binary variable $InOut_{ab}$ to determine the order.
\end{itemize}

\subsubsection{Blocking Constraints}
\begin{figure}[htbp]
    \centering
    \includegraphics[width=0.6\textwidth]{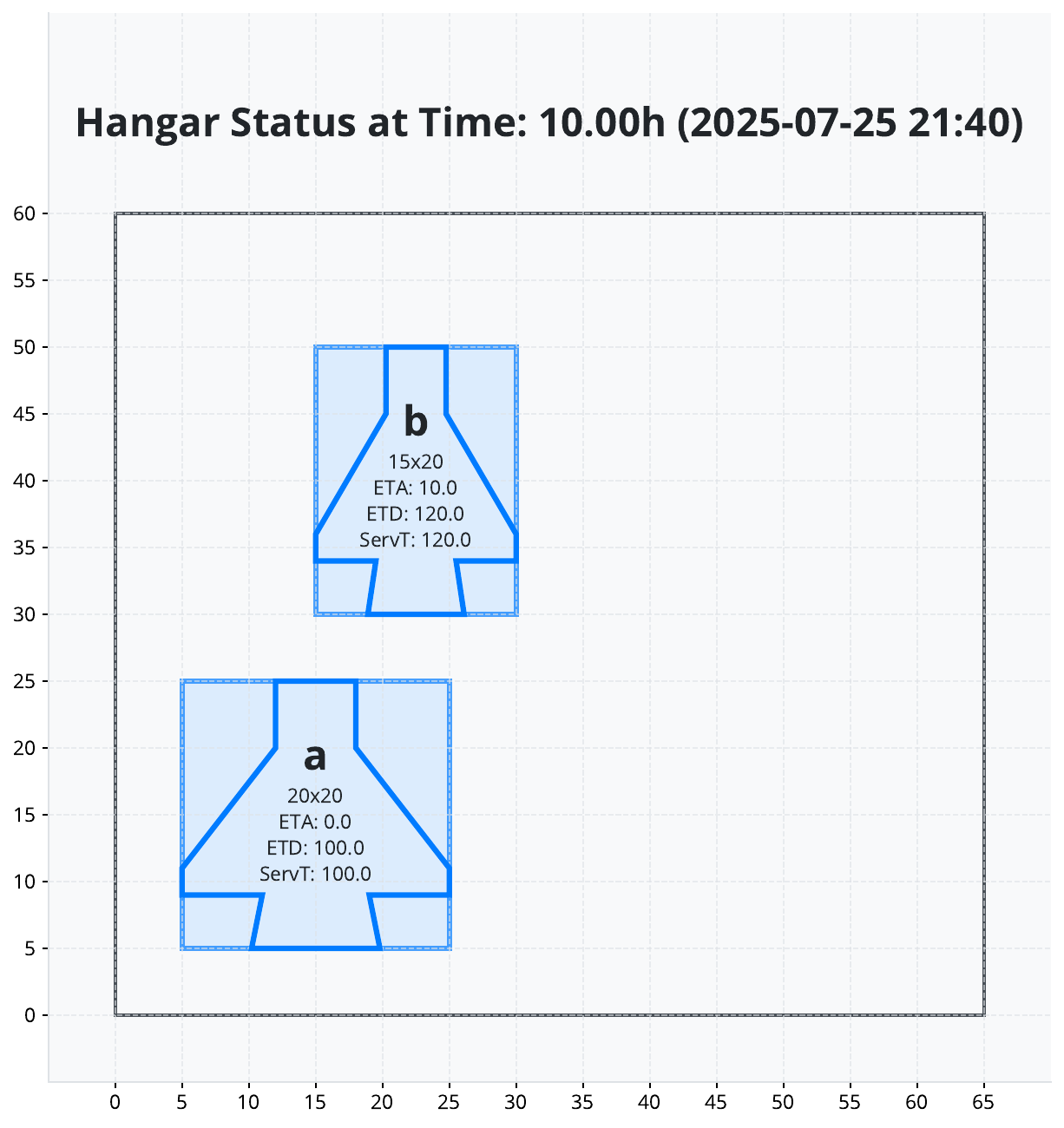} 
    \caption{An illustration of a blocking scenario. Aircraft `b` is parked at a higher Y-coordinate, obstructing aircraft `a` from reaching the hangar exit (assumed at the top). For the departure case, this situation arises if `a` arrived first and `b` subsequently parked in the blocking position; consequently, `b` must depart before `a` can exit. For the arrival case, if `b` is already present, `a` cannot move into its depicted position until `b` has departed.}
    \label{fig:blocking}
\end{figure}

These constraints model simplified logistical pathways to prevent deadlocks where one aircraft's movement is obstructed by another. A simplified blocking scenario is illustrated in Figure~\ref{fig:blocking}.
\begin{align}
& \begin{aligned}
& Roll\_out_a \ge Roll\_out_b + \varepsilon_t \\
& - M_T \cdot ((1 - Above_{ba}) + Right_{ab} + Right_{ba} \\
& + (1 - InIn_{ab})),
\end{aligned} && \parbox{3.5cm}{\raggedright$\forall a, b \in A, a \neq b$} \label{eq:block1} \\
& \begin{aligned}
& Roll\_in_a \ge Roll\_out_b + \varepsilon_t \\
& - M_T \cdot ((1 - Above_{ba}) + Right_{ab} + Right_{ba} + InIn_{ab}),
\end{aligned} \label{eq:block2} \\
& \qquad \qquad \qquad \qquad \qquad \forall a, b \in A, a \neq b, \text{ where } a \in F \lor b \in F \nonumber
\end{align}
\textbf{Explanation:}
\begin{itemize}
    \item \textbf{Constraint \eqref{eq:block1} (Departure Blocking):} This constraint prevents a departure deadlock. Consider the scenario in Figure~\ref{fig:blocking}, where aircraft `b` is positioned at a higher Y-coordinate than aircraft `a`, directly obstructing its path to the hangar exit. This situation is captured mathematically when `b` is above `a` ($Above_{ba}=1$), they are in the same vertical lane (neither is to the right or left of the other, so $Right_{ab}=0$ and $Right_{ba}=0$), and `a` arrived before `b` ($InIn_{ab}=1$). Under these specific conditions, the Big-M term in the constraint becomes zero, activating the rule $Roll\_out_a \ge Roll\_out_b + \varepsilon_t$. Using the example data from the figure, `b`'s departure ($Roll\_out_b$) is scheduled at time 120. This rule forces `a`'s departure to occur after 120, creating a departure delay ($D^{Dep}_a$) since its original $ETD_a$ was 100. The model must then weigh the cost of this delay against other alternatives, such as placing `b` in a non-blocking position.
    \item \textbf{Constraint \eqref{eq:block2} (Arrival Blocking):} This constraint handles a similar scenario for arrivals. Imagine aircraft `b` is already parked as shown in Figure~\ref{fig:blocking}. If the model attempts to schedule a new aircraft `a` to arrive and park in the position shown (below `b`), the same blocking condition exists ($Above_{ba}=1$, same lane). The constraint activates, enforcing $Roll\_in_a \ge Roll\_out_b + \varepsilon_t$. This means `a` is not allowed to enter the hangar until the blocking aircraft `b` has departed, potentially causing a significant arrival delay ($D^{Arr}_a$).
\end{itemize}

\subsection{Initial Conditions and Variable Domains}

\subsubsection{Fixed Initial Values}
These constraints initialize the model based on the state of the hangar at the beginning of the planning horizon.
\begin{align}
& Accept_c = 1, && \forall c \in C \label{eq:init_accept} \\
& X_c = X^{init}_c, && \forall c \in C \label{eq:init_pos_x} \\
& Y_c = Y^{init}_c, && \forall c \in C \label{eq:init_pos_y} \\
& Roll\_in_c = 0, && \forall c \in C \label{eq:init_rollin} \\
& InIn_{cf} = 1, && \forall c \in C, f \in F \label{eq:init_inin_cf} \\
& InIn_{fc} = 0, && \forall f \in F, c \in C \label{eq:init_inin_fc} \\
& InIn_{cd} = 1, && \forall c, d \in C, c \neq d \label{eq:init_inin_cd}
\end{align}

\subsubsection{Variable Type Definitions}
\begin{align}
& X_a, Y_a, Roll\_in_a, Roll\_out_a, D^{Arr}_f, D^{Dep}_a \ge 0 \label{eq:domain_pos} \\
& Accept_a, Right_{ab}, Above_{ab}, \dots \in \{0, 1\} \label{eq:domain_bin}
\end{align}
\textbf{Explanation:}
\begin{itemize}
    \item \textbf{Constraints \eqref{eq:init_accept}-\eqref{eq:init_rollin}:} These equations fix the status of aircraft already in the hangar: they are considered accepted, their positions are known, and their roll-in time is set to zero relative to the start of the plan.
    \item \textbf{Constraints \eqref{eq:init_inin_cf}-\eqref{eq:init_inin_cd}:} These constraints pre-define the roll-in sequence to be consistent with reality: all current aircraft ($c$) are considered to have rolled in before any future aircraft ($f$). They also establish a fixed, albeit arbitrary, roll-in order among the current aircraft to satisfy sequencing logic.
    \item \textbf{Constraints \eqref{eq:domain_pos} and \eqref{eq:domain_bin}:} These define the domains for the decision variables, ensuring continuous variables are non-negative and binary variables take values of 0 or 1.
\end{itemize}

\section{Heuristic Solution Approach}
While the MILP model guarantees optimality, its computational requirements can become significant for very large-scale instances or in time-critical operational contexts. To provide a performance benchmark and a rapid decision-making alternative, we developed an \textbf{Automated Constructive Heuristic (ACH)}. The purpose of the ACH is not to compete with sophisticated metaheuristics (e.g., genetic algorithms or simulated annealing), but rather to serve as a robust baseline that mimics a logical, priority-driven human planning process. By comparing the MILP's optimal solutions against the outcomes of this fast, greedy approach, we can rigorously quantify the economic value of achieving true optimality.

\subsection{Algorithm Design and Rationale}
The ACH is a deterministic, single-pass constructive heuristic implemented in Python. It makes sequential, irrevocable decisions, building a feasible solution by processing one aircraft at a time. Its design is grounded in a set of logical rules that ensure every decision respects the core operational constraints defined in the mathematical model. This direct correspondence ensures that the heuristic produces physically and logistically viable schedules that are directly comparable to the MILP solutions. The algorithm operates in two primary phases: prioritization and sequential placement.

\subsubsection{Phase 1: Prioritization Rule}
The cornerstone of the heuristic is its prioritization logic, which determines the sequence in which new aircraft requests ($f \in F$) are considered for placement. A well-defined priority queue is essential for a greedy algorithm, as the initial decisions heavily constrain subsequent options. All new aircraft are sorted into a priority list, $L$, based on the following lexicographical criteria, applied in order:
\begin{enumerate}
    \item \textbf{Rejection Penalty ($P^{Rej}_f$):} Descending order. Aircraft with higher rejection penalties are prioritized to minimize the most significant potential costs first. This reflects the high economic or strategic value of certain maintenance tasks.
    \item \textbf{Expected Time of Arrival ($ETA_f$):} Ascending order. Among aircraft with equal rejection penalties, those scheduled to arrive earlier are considered first. This rule aims to adhere to the original schedule as closely as possible.
    \item \textbf{Service Time ($ServT_f$):} Ascending order. As a final tie-breaker, aircraft requiring less time in the hangar are prioritized. This heuristic choice aims to increase hangar throughput by processing quicker tasks first, potentially freeing up space for subsequent, longer-duration jobs.
\end{enumerate}
This multi-attribute rule creates a deterministic and operationally sensible sequence that guides the heuristic's greedy decisions.

\subsubsection{Phase 2: Sequential Placement and Scheduling}
After establishing the priority list $L$, the algorithm iterates through it, attempting to find a feasible spatiotemporal slot for each aircraft. For a given aircraft, the heuristic initiates a search for a valid placement starting at its earliest possible roll-in time, which is its $ETA_f$. The search space is explored in two dimensions: time and space.

The core of this phase is the feasibility check for a potential placement. For a given coordinate pair $(x, y)$ and a roll-in time $t$, this function rigorously validates the decision against the set of fixed, previously scheduled aircraft. It enforces the same critical constraints as the MILP model:
\begin{itemize}
    \item \textbf{Hangar Boundaries:} The aircraft's footprint, including the safety `Buffer`, must be entirely within the hangar dimensions ($HW, HL$).
    \item \textbf{Spatial Non-overlapping:} The new aircraft's buffered bounding box cannot overlap with that of any other aircraft already scheduled to be in the hangar during the same time interval.
    \item \textbf{Temporal Separation:} All roll-in and roll-out events must be separated by the minimum time gap, $\varepsilon_t$.
    \item \textbf{Blocking Rules:} The heuristic enforces path-blocking logic. An aircraft cannot be placed in a position that blocks the exit path of an already-parked aircraft, nor can it be scheduled to arrive in a position whose entrance path is blocked by an existing aircraft until that aircraft has departed.
\end{itemize}

If, at a given time $t$, one or more valid placements are found, the algorithm does not merely select the first one. Instead, it performs an exhaustive spatial scan of the entire hangar grid, collects all feasible positions into a candidate set, and then selects the best spot among them. The selection criterion is the minimization of the sum of coordinates $(x+y)$, a greedy choice that promotes a compact and organized layout by favoring positions closer to the hangar's origin.

If no valid placement is found at time $t$, the algorithm incrementally delays the potential roll-in time by a small step, $\varepsilon_t$, and re-evaluates the entire hangar space at the new time. This process continues until a feasible slot is identified or the potential arrival delay becomes economically unviable.

The decision to reject an aircraft is guided by a straightforward economic rationale: the search continues only as long as the cumulative cost of delaying its arrival is lower than the fixed penalty for rejection. To formalize this, the algorithm determines the economic break-even point, where the total delay cost equals the rejection penalty. This yields a threshold on the maximum allowable arrival delay, $D^{Arr}_{max}$, computed as follows:
\begin{equation*}
D^{Arr}_{max} = \frac{P^{Rej}_f}{P^{Arr}_f}
\end{equation*}
Accordingly, the latest admissible arrival time, $t_{max}$, is the expected time of arrival ($ETA_f$) plus this maximum delay:
\begin{equation*}
t_{max} = ETA_f + \frac{P^{Rej}_f}{P^{Arr}_f}
\end{equation*}

\begin{figure}[htbp!]
    \centering
    \includegraphics[width=\textwidth]{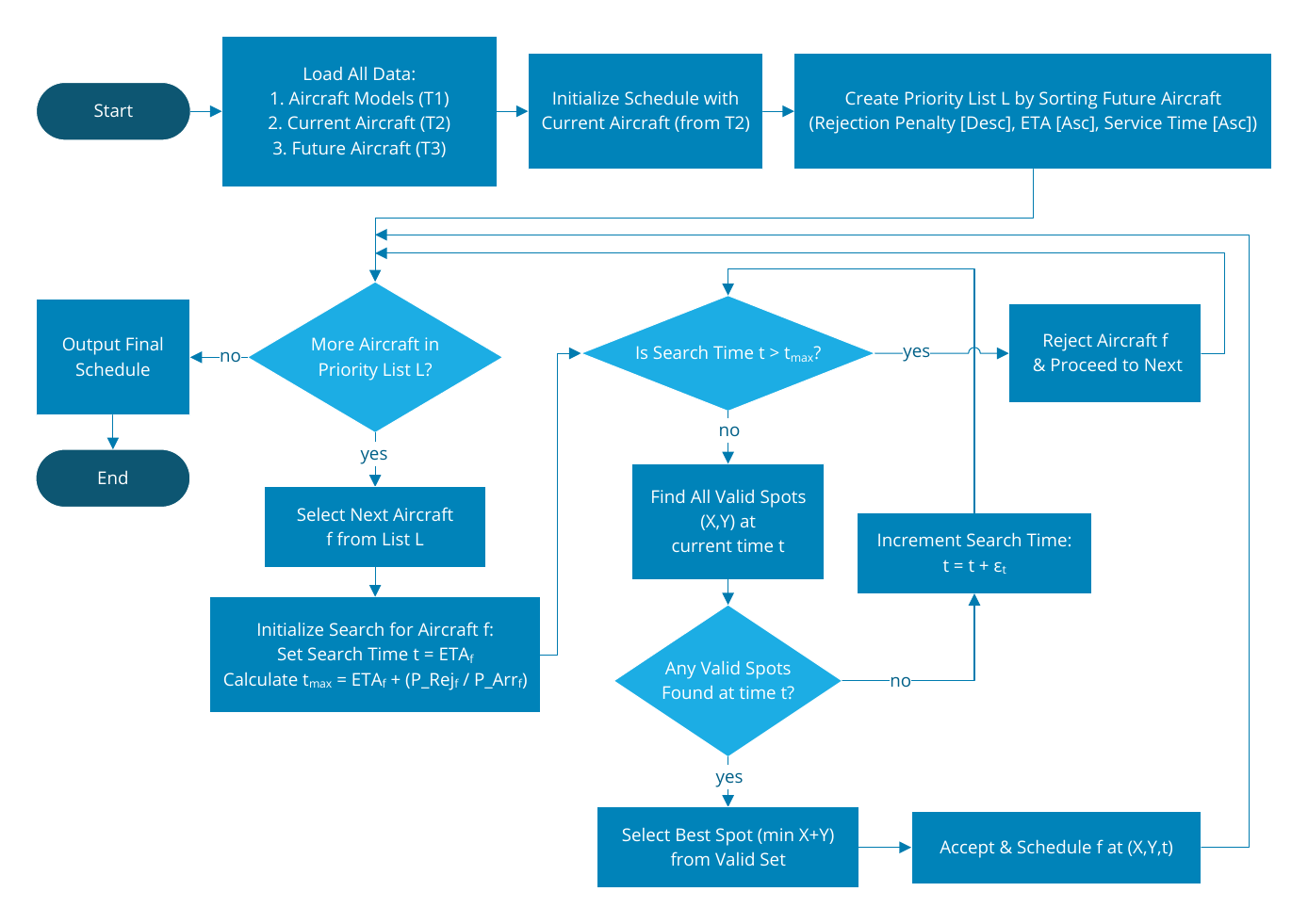}
    \caption{Flowchart of the Automated Constructive Heuristic (ACH). The "Valid Spot" check encapsulates all spatial, temporal, and blocking constraints, including hangar boundaries, non-overlapping rules, minimum time gaps, and path blocking logic.}
    \label{fig:heuristic_flowchart}
\end{figure}

The search for a valid spot is terminated and the aircraft is rejected if the search time $t$ exceeds this calculated $t_{max}$. This mechanism provides a smart, cost-aware trade-off between delaying and rejecting service. The heuristic's logic is visualized in Figure~\ref{fig:heuristic_flowchart}.

\section{Computational Study and Results}
This section presents a comprehensive computational study designed to rigorously evaluate the performance of our proposed continuous-time MILP model. Our evaluation is structured as a multi-stage investigation with a clear narrative arc. We first establish our model's superiority by benchmarking it against the current state-of-the-art. Building on this foundation, we conduct an extensive scalability analysis to map the model's performance boundaries across a wide range of systematically generated instances. Finally, we perform a focused sensitivity analysis to probe the model's resilience to temporal congestion, a known bottleneck for alternative formulations. This structured approach allows us to not only validate our model but also to derive actionable insights into its practical applicability.

\subsection{Experimental Design and Testbed Generation}
All experiments were conducted on a machine with an Intel Core i7-10750H CPU @ 2.60GHz and 16 GB of RAM. The MILP model was implemented in GAMS and solved using the CPLEX solver. The Automated Constructive Heuristic (ACH) was implemented in Python.

\subsubsection{Instance Generation Framework for Scalability and Sensitivity Tests}
To create a robust and representative testbed, we developed a systematic framework for generating the \texttt{RND}, \texttt{INC}, and \texttt{CON} instance sets. These instances were configured with a set of default parameters reflecting a plausible operational scenario. The hangar dimensions were set to a width (HW) of 65 meters and a length (HL) of 60 meters, with a 5-meter safety buffer, representing a moderately sized MRO facility. The minimum time gap between movements ($\varepsilon_t$) was set to 0.1 hours (6 minutes). All generated instances also begin with an initial state of two pre-existing aircraft in the hangar to simulate ongoing operations.

A key feature of our generation methodology is the incorporation of realistic, correlated heterogeneity:
\begin{enumerate}
    \item \textbf{Controlled Scalability:} The time horizon for new aircraft arrivals ($ETA_{max}$) scales linearly with the number of requests (N). This approach preserves the relative event density across all problem sizes, ensuring that our scalability analysis measures the pure effect of increasing N, without being confounded by changes in temporal congestion.
    
    \item \textbf{Correlated Economic and Physical Attributes:} To create more realistic decision-making challenges, economic attributes are correlated with physical ones. Each aircraft is assigned a type (M\_ID from 1 to 8, representing size). The probability of an aircraft being designated a 'VIP' is then determined by its type, with larger aircraft having a higher likelihood of being VIPs. This reflects the real-world scenario where larger assets often have higher operational priority. However, the framework still allows for exceptions, such as large non-VIP aircraft or small, high-priority VIP jets, enriching the complexity of the trade-offs.
    
    \item \textbf{Heterogeneous Penalties:} Penalty costs are not fixed. They are randomized within specified ranges that differ for VIP and non-VIP aircraft. Furthermore, the base rejection penalty is also correlated with the aircraft's type, adding another layer of realistic complexity that moves beyond simple, uniform cost structures.
\end{enumerate}

\subsubsection{Experimental Campaign Outline}
Our computational campaign was divided into distinct experiments to test different facets of the model's performance:
\begin{itemize}
    \item \textbf{Benchmark Comparison:} We test our model on the \texttt{Case2015} instances from \citet{qin2019} to provide a direct comparison against the event-based state-of-the-art.
    
    \item \textbf{Scalability Analysis:} To assess the model's scalability, we performed extensive tests on two distinct instance sets: \texttt{RND} for evaluating performance on random structural variations, and \texttt{INC} for measuring the impact of controlled, incremental growth.
    For the \texttt{RND} set, small-to-medium scale instances were generated for N from 5 to 40 (in steps of 5), while large-scale instances covered N from 60 to 160 (in steps of 20). To ensure statistical robustness, three independent replications were created for each size.
    For the \texttt{INC} set, small-to-medium instances were generated in finer steps, with N ranging from 6 to 40 (in steps of 2), while large-scale instances mirrored the \texttt{RND} set, covering N from 60 to 160 (in steps of 20).
    In all scalability tests, small-to-medium instances were solved to proven optimality, whereas large-scale instances were run with a 3600-second time limit to evaluate performance under practical computational constraints.

    \item \textbf{Sensitivity Analysis:} To measure performance degradation under systematically increasing levels of temporal congestion, we conducted a focused analysis. For this, we used three distinct base instances of size N=20, and for each one, we generated several scenarios by systematically varying the congestion level.
\end{itemize}

\subsection{Direct Comparison with State-of-the-Art Benchmark}
To ground our contributions within the existing literature, we begin our analysis with a direct comparison against the event-based model of \citet{qin2019}. We use their \texttt{Case2015} instances, which are based on a real-world case study and serve as a recognized benchmark. These instances utilize a different parameter set, featuring a larger hangar and fixed penalty costs, providing a robust test of our model's architectural advantages. The results of this head-to-head comparison are presented in Table~\ref{tab:case2015_results}.

\begin{table}[htbp!]
\centering
\caption{Performance comparison on the Case2015 benchmark instances from \citet{qin2019}.}
\label{tab:case2015_results}
\begin{tabular}{l r r r l r l}
\toprule
\multicolumn{4}{c}{\thead{Our Model \\ (Continuous-Time)}} & \multicolumn{3}{c}{\thead{Qin et al., 2019 \\ (Event-Based)}} \\
\cmidrule(lr){1-4} \cmidrule(lr){5-7}
\thead{Instance \\ Name} & \thead{Total \\ Time (s)} & \thead{Total \\ Cost} & \thead{Final \\ Gap (\%)} & \thead{Instance \\ Name} & \thead{Reported \\ Time (s)} & \thead{Reported \\ Gap (\%)} \\
\midrule
Case15-C9 & 0.11 & 160 & 0 & 3\_2\_9 & 3600   & 100 (Failed) \\
Case15-S9 & 0.09 & 320 & 0 & 4\_3\_9 & 457.86 & 0 \\
Case15-E8 & 0.28 & 160 & 0 & 4\_1\_8 & 2.00   & 0 \\
\bottomrule
\end{tabular}
\end{table}

\begin{figure}[htbp!]
    \centering
    \includegraphics[width=0.8\textwidth]{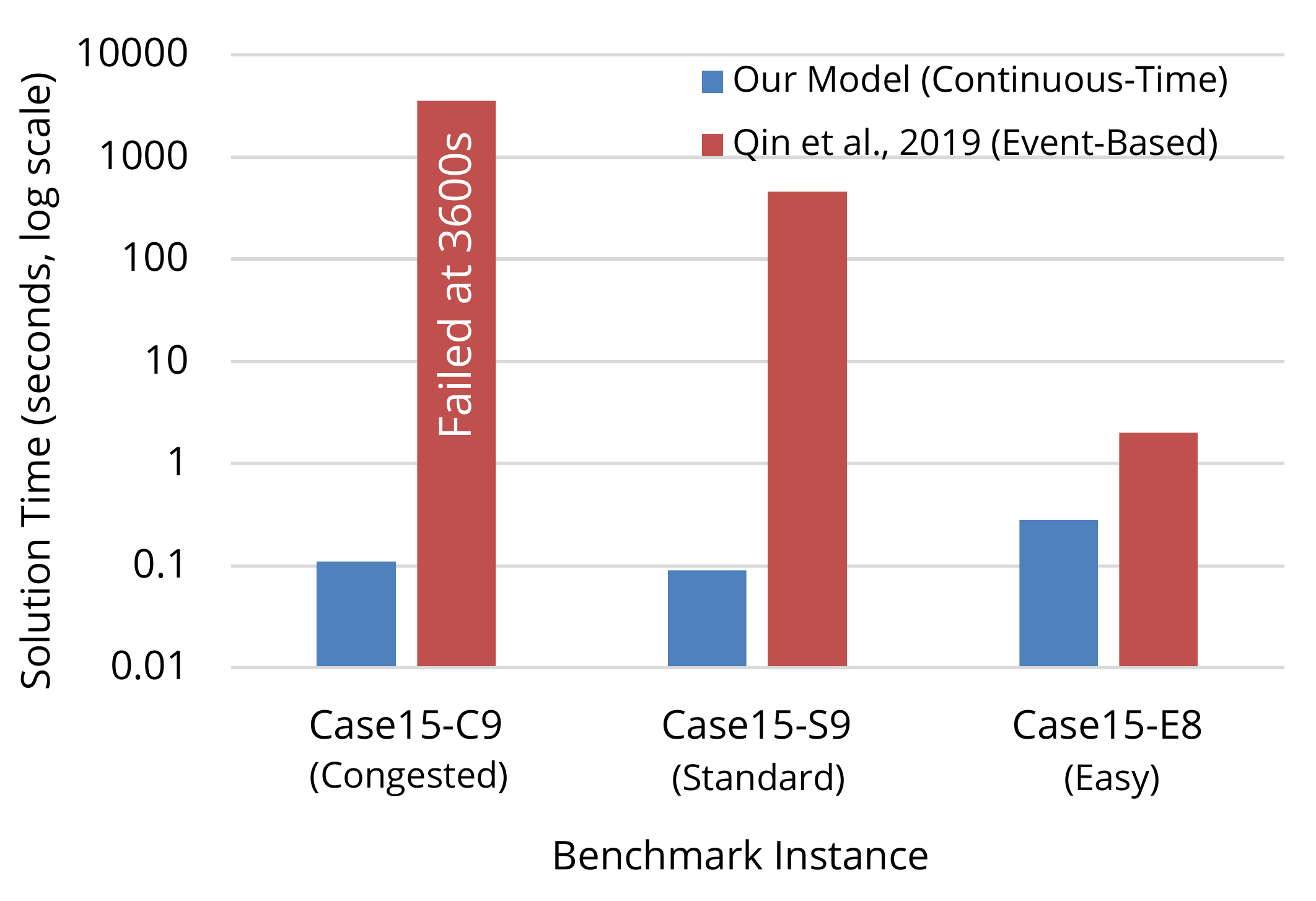}
    \caption{Solution time comparison on benchmark instances from \citet{qin2019}. The y-axis is on a logarithmic scale to visualize the orders-of-magnitude performance difference.}
    \label{fig:benchmark_comparison}
\end{figure}

The results in Table~\ref{tab:case2015_results}, visually summarized in Figure~\ref{fig:benchmark_comparison}, demonstrate a transformative improvement in computational efficiency. The most critical finding is for instance \texttt{Case15-C9}, a congested scenario that the event-based model failed to solve within a one-hour time limit. Our continuous-time model solves this previously intractable instance to proven optimality in a mere \textbf{0.11 seconds}. This outcome is not an anomaly; for the standard (\texttt{Case15-S9}) and easy (\texttt{Case15-E8}) instances, our model is over 5,000 and 7 times faster, respectively. This comprehensive outperformance across all problem types provides definitive evidence that our continuous-time formulation represents a fundamentally more powerful and efficient approach, motivating a deeper exploration of its capabilities on a broader and more complex testbed.

\subsection{Scalability and Performance Analysis}
Having established the model's fundamental superiority on established benchmarks, we now turn to a more rigorous and extensive evaluation of its performance limits. This analysis is crucial for understanding the model's practical applicability in real-world settings, which often involve a larger number of aircraft and more complex operational constraints. We investigate scalability through two lenses: performance on randomly generated instances with diverse characteristics, and performance on instances that grow incrementally in size.

To provide context for the performance results, Table~\ref{tab:milp_specs_rnd} and Table~\ref{tab:milp_specs_inc} first detail the structural growth of the MILP model itself for both the random (\texttt{RND}) and incremental (\texttt{INC}) instance sets, respectively. These tables illustrate how the number of equations, variables, and discrete variables increases polynomially with the number of aircraft, highlighting the underlying computational challenge.

\begin{table}[htbp!]
\centering
\sisetup{group-separator={,}}
\caption{MILP model specifications for the generated RND instances, showing the growth in model size.}
\label{tab:milp_specs_rnd}
\begin{tabular}{l S[table-format=3.0] S[table-format=6.0] S[table-format=6.0] S[table-format=6.0]}
\toprule
\thead{Instance Name} & {\thead{Total \\ Aircraft}} & {\thead{Equations}} & {\thead{Variables}} & {\thead{Discrete \\ Variables}} \\
\midrule
RND-N005  & 7   & 441     & 277     & 212     \\
RND-N010  & 12  & 1391    & 807     & 692     \\
RND-N015  & 17  & 2866    & 1612    & 1447    \\
RND-N020  & 22  & 4866    & 2692    & 2477    \\
RND-N025  & 27  & 7391    & 4047    & 3782    \\
RND-N030  & 32  & 10441   & 5677    & 5362    \\
RND-N035  & 37  & 14016   & 7582    & 7217    \\
RND-N040  & 42  & 18116   & 9762    & 9347    \\
RND-N060  & 62  & 39766   & 21232   & 20617   \\
RND-N080  & 82  & 69816   & 37102   & 36287   \\
RND-N100  & 102 & 108266  & 57372   & 56357   \\
RND-N120  & 122 & 155116  & 82042   & 80827   \\
RND-N140  & 142 & 210366  & 111112  & 109697  \\
RND-N160  & 162 & 274016  & 144582  & 142967  \\
\bottomrule
\end{tabular}
\end{table}

\begin{table}[htbp!]
\centering
\sisetup{group-separator={,}, group-minimum-digits=4}
\caption{MILP model specifications for the generated INC instances, showing the controlled growth in model size.}
\label{tab:milp_specs_inc}
\begin{tabular}{l S[table-format=3.0] S[table-format=6.0] S[table-format=6.0] S[table-format=6.0]}
\toprule
\thead{Instance Name} & {\thead{Total \\ Aircraft}} & {\thead{Equations}} & {\thead{Variables}} & {\thead{Discrete \\ Variables}} \\
\midrule
INC-N006  & 8   & 589     & 361     & 286     \\
INC-N008  & 10  & 948     & 562     & 467     \\
INC-N010  & 12  & 1391    & 807     & 692     \\
INC-N012  & 14  & 1918    & 1096    & 961     \\
INC-N014  & 16  & 2529    & 1429    & 1274    \\
INC-N016  & 18  & 3224    & 1806    & 1631    \\
INC-N018  & 20  & 4003    & 2227    & 2032    \\
INC-N020  & 22  & 4866    & 2692    & 2477    \\
INC-N022  & 24  & 5813    & 3201    & 2966    \\
INC-N024  & 26  & 6844    & 3754    & 3499    \\
INC-N026  & 28  & 7959    & 4351    & 4076    \\
INC-N028  & 30  & 9158    & 4992    & 4697    \\
INC-N030  & 32  & 10441   & 5677    & 5362    \\
INC-N032  & 34  & 11808   & 6406    & 6071    \\
INC-N034  & 36  & 13259   & 7179    & 6824    \\
INC-N036  & 38  & 14794   & 7996    & 7621    \\
INC-N038  & 40  & 16413   & 8857    & 8462    \\
INC-N040  & 42  & 18116   & 9762    & 9347    \\
INC-N060  & 62  & 39766   & 21232   & 20617   \\
INC-N080  & 82  & 69816   & 37102   & 36287   \\
INC-N100  & 102 & 108266  & 57372   & 56357   \\
INC-N120  & 122 & 155116  & 82042   & 80827   \\
INC-N140  & 142 & 210366  & 111112  & 109697  \\
INC-N160  & 162 & 274016  & 144582  & 142967  \\
\bottomrule
\end{tabular}
\end{table}

\subsubsection{Performance on Small-to-Medium Scale Instances (N \texorpdfstring{$\leq$}{<=} 40)}
We first evaluate the model on small-to-medium instances, where we seek to find and prove optimal solutions. Table~\ref{tab:small_scale_random_results} and Table~\ref{tab:small_scale_inc_results} present detailed performance comparisons for the \texttt{RND} and \texttt{INC} sets, respectively. These tables not only report the optimal solution quality and time but also include intermediate results for reaching near-optimality (10\% and 5\% gaps), offering insights into the speed-quality trade-off.

\begin{figure}[htbp!]
    \centering
    \includegraphics[width=0.8\textwidth]{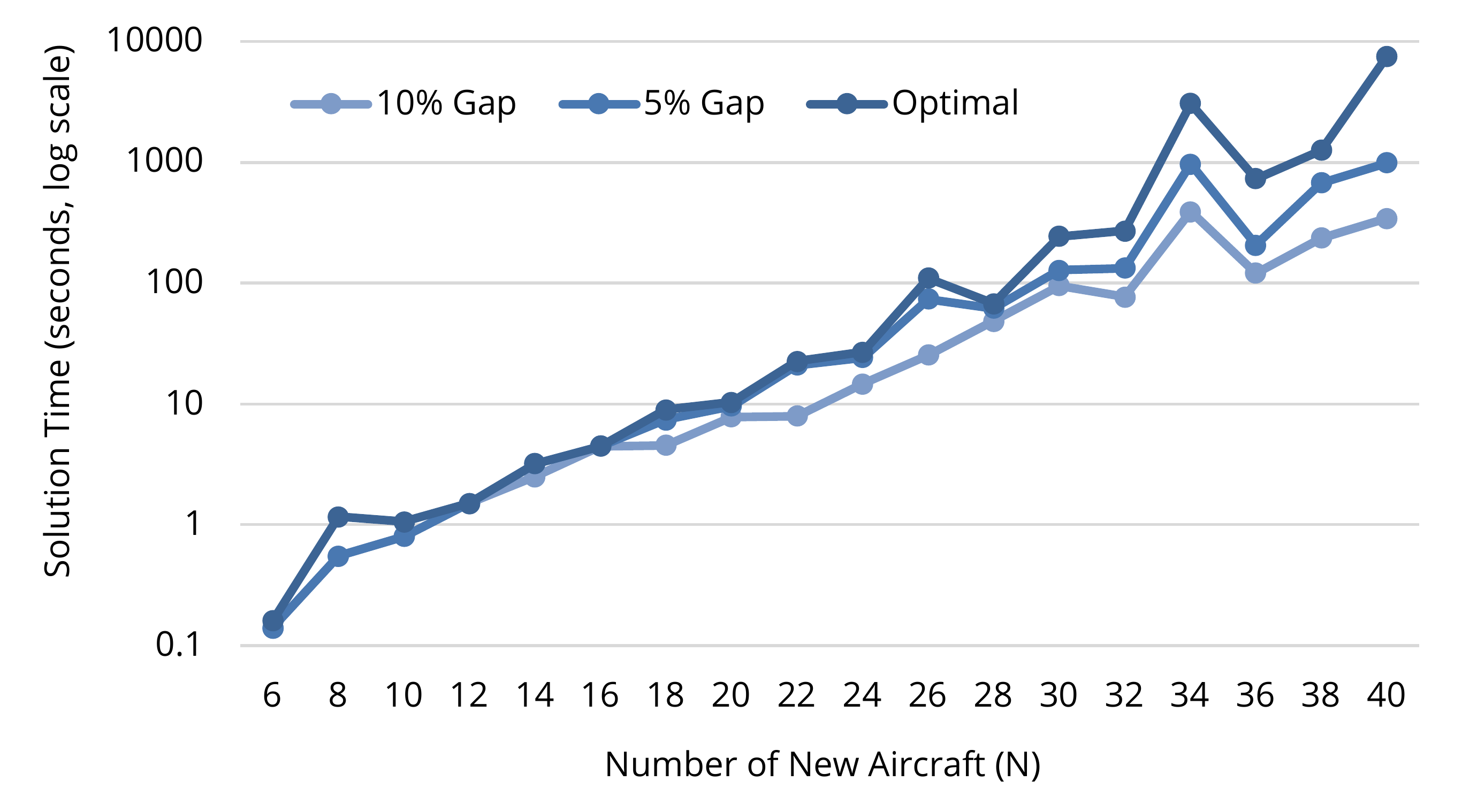} 
    \caption{Solution time required to reach different optimality gaps for the incremental (\texttt{INC}) instances. The y-axis is on a logarithmic scale.}
    \label{fig:scalability_gaps}
\end{figure}

\begin{sidewaystable}[htbp!]
\centering
\scriptsize 
\caption{Detailed performance comparison on small-to-medium scale random instances (N = 5 to 40).}
\label{tab:small_scale_random_results}
\sisetup{group-separator={,}, group-minimum-digits=4}
\begin{tabular}{l *{10}{S[table-format=5.0]}}
\toprule
\multirow{2}{*}{\thead{Instance Name}} & \multicolumn{2}{c}{\thead{MILP (10\% Gap)}} & \multicolumn{2}{c}{\thead{MILP (5\% Gap)}} & \multicolumn{2}{c}{\thead{Optimal MILP}} & \multicolumn{2}{c}{\thead{ACH Heuristic}} & \multicolumn{2}{c}{\thead{Accepted Aircraft}} \\
\cmidrule(lr){2-3} \cmidrule(lr){4-5} \cmidrule(lr){6-7} \cmidrule(lr){8-9} \cmidrule(lr){10-11}
& {\thead{Time (s)}} & {\thead{Total \\ Cost}} & {\thead{Time (s)}} & {\thead{Total \\ Cost}} & {\thead{Time (s)}} & {\thead{Total \\ Cost}} & {\thead{Total \\ Cost}} & {\thead{Time (s)}} & {\thead{ACH}} & {\thead{MILP}} \\
\midrule
RND-N005-I01 & 0    & 4791   & 0    & 4791   & 0      & 4791   & 4791   & 6  & 3 & 3 \\
RND-N005-I02 & 0    & 3568   & 0    & 3568   & 0      & 3568   & 3777   & 5  & 4 & 5 \\
RND-N005-I03 & 0    & 11876  & 0    & 11876  & 0      & 11876  & 16816  & 8  & 3 & 4 \\
\midrule
RND-N010-I01 & 1    & 6730   & 1    & 6730   & 1      & 6730   & 7912   & 11 & 5 & 6 \\
RND-N010-I02 & 1    & 10902  & 1    & 10902  & 1      & 10902  & 28739  & 15 & 6 & 8 \\
RND-N010-I03 & 1    & 12784  & 1    & 12784  & 1      & 12784  & 12784  & 9  & 6 & 6 \\
\midrule
RND-N015-I01 & 3    & 28062  & 3    & 27779  & 3      & 27185  & 49079  & 20 & 7 & 10 \\
RND-N015-I02 & 3    & 20093  & 8    & 19721  & 12     & 19721  & 32674  & 20 & 7 & 12 \\
RND-N015-I03 & 2    & 26130  & 2    & 26130  & 2      & 26130  & 36173  & 18 & 6 & 9 \\
\midrule
RND-N020-I01 & 9    & 18970  & 9    & 18873  & 14     & 18873  & 34308  & 21 & 11 & 15 \\
RND-N020-I02 & 4    & 28559  & 5    & 27667  & 8      & 27667  & 35497  & 24 & 9  & 13 \\
RND-N020-I03 & 3    & 33921  & 4    & 33776  & 5      & 33771  & 44436  & 24 & 10 & 15 \\
\midrule
RND-N025-I01 & 37   & 22975  & 45   & 22569  & 50     & 22569  & 39366  & 27 & 12 & 16 \\
RND-N025-I02 & 15   & 29687  & 16   & 29479  & 20     & 29479  & 41695  & 27 & 11 & 16 \\
RND-N025-I03 & 6    & 47179  & 7    & 46132  & 7      & 46132  & 87276  & 35 & 10 & 14 \\
\midrule
RND-N030-I01 & 31   & 33116  & 34   & 32501  & 41     & 32253  & 41551  & 32 & 12 & 16 \\
RND-N030-I02 & 38   & 60444  & 39   & 57531  & 80     & 56103  & 90928  & 39 & 12 & 19 \\
RND-N030-I03 & 36   & 34002  & 70   & 32879  & 170    & 32879  & 46696  & 30 & 15 & 20 \\
\midrule
RND-N035-I01 & 163  & 43877  & 284  & 43877  & 526    & 43854  & 67215  & 36 & 18 & 24 \\
RND-N035-I02 & 67   & 64753  & 91   & 63991  & 105    & 61855  & 89069  & 43 & 12 & 20 \\
RND-N035-I03 & 51   & 41164  & 52   & 39409  & 85     & 39092  & 63828  & 35 & 18 & 23 \\
\midrule
RND-N040-I01 & 172  & 75864  & 1037 & 74456  & 58837  & 74456  & 121782 & 47 & 17 & 23 \\
RND-N040-I02 & 2333 & 42367  & 8994 & 39712  & 31709  & 41804  & 62121  & 41 & 19 & 24 \\
RND-N040-I03 & 89   & 65047  & 143  & 63981  & 520    & 63746  & 96065  & 50 & 14 & 22 \\
\bottomrule
\end{tabular}
\end{sidewaystable}

\begin{sidewaystable}[htbp!]
\centering
\footnotesize 
\caption{Detailed performance comparison on small-to-medium scale incremental instances (N = 6 to 40).}
\label{tab:small_scale_inc_results}
\sisetup{group-separator={,}, group-minimum-digits=4}
\begin{tabular}{l *{10}{S[table-format=5.0]}}
\toprule
\multirow{2}{*}{\thead{Instance Name}} & \multicolumn{2}{c}{\thead{MILP (10\% Gap)}} & \multicolumn{2}{c}{\thead{MILP (5\% Gap)}} & \multicolumn{2}{c}{\thead{Optimal MILP}} & \multicolumn{2}{c}{\thead{ACH Heuristic}} & \multicolumn{2}{c}{\thead{Accepted Aircraft}} \\
\cmidrule(lr){2-3} \cmidrule(lr){4-5} \cmidrule(lr){6-7} \cmidrule(lr){8-9} \cmidrule(lr){10-11}
& {\thead{Time (s)}} & {\thead{Total \\ Cost}} & {\thead{Time (s)}} & {\thead{Total \\ Cost}} & {\thead{Total \\ Time (s)}} & {\thead{Total \\ Cost}} & {\thead{Total \\ Cost}} & {\thead{Time (s)}} & {\thead{ACH}} & {\thead{MILP}} \\
\midrule
INC-N006 & 0   & 11659 & 0   & 11659 & 0    & 11659 & 29612  & 9  & 5  & 6  \\
INC-N008 & 1   & 13215 & 1   & 13215 & 1    & 13215 & 25637  & 11 & 5  & 7  \\
INC-N010 & 1   & 15552 & 1   & 15552 & 1    & 15552 & 27974  & 13 & 5  & 7  \\
INC-N012 & 2   & 22568 & 2   & 22568 & 2    & 22568 & 38222  & 18 & 5  & 8  \\
INC-N014 & 3   & 26217 & 3   & 26198 & 3    & 26198 & 49351  & 20 & 7  & 8  \\
INC-N016 & 4   & 31698 & 4   & 31698 & 4    & 31698 & 50147  & 22 & 6  & 9  \\
INC-N018 & 5   & 37457 & 7   & 37388 & 9    & 37198 & 56626  & 26 & 6  & 10 \\
INC-N020 & 8   & 38499 & 10  & 38499 & 10   & 38309 & 71047  & 31 & 7  & 11 \\
INC-N022 & 8   & 38788 & 21  & 38598 & 23   & 38593 & 72002  & 32 & 8  & 13 \\
INC-N024 & 15  & 39756 & 24  & 39756 & 27   & 39749 & 65062  & 33 & 8  & 14 \\
INC-N026 & 25  & 45720 & 74  & 45525 & 110  & 45525 & 82116  & 36 & 10 & 14 \\
INC-N028 & 48  & 45939 & 62  & 45749 & 68   & 45749 & 84117  & 39 & 10 & 16 \\
INC-N030 & 96  & 48076 & 128 & 48069 & 245  & 48069 & 86490  & 41 & 10 & 16 \\
INC-N032 & 77  & 49246 & 133 & 49246 & 272  & 49241 & 88657  & 45 & 10 & 17 \\
INC-N034 & 390 & 49543 & 971 & 49543 & 3073 & 49543 & 89486  & 46 & 11 & 19 \\
INC-N036 & 121 & 56259 & 207 & 55050 & 734  & 55043 & 103949 & 49 & 13 & 20 \\
INC-N038 & 237 & 61898 & 677 & 61898 & 1261 & 61371 & 109632 & 53 & 12 & 21 \\
INC-N040 & 342 & 66702 & 992 & 66246 & 7556 & 66246 & 115755 & 56 & 12 & 22 \\
\bottomrule
\end{tabular}
\end{sidewaystable}

A key insight emerges from these tables: the value of optimality extends beyond mere cost reduction. By navigating complex spatio-temporal trade-offs with a global perspective, the MILP model consistently accepts a greater number of aircraft than the myopic heuristic. For instance, in \texttt{RND-N030-I03}, the optimal plan accommodates 20 aircraft versus the heuristic's 15. This demonstrates that the model's primary contribution is unlocking significant hidden operational capacity, a critical factor for MRO profitability. Furthermore, the results highlight the model's strategic flexibility. As illustrated in Figure~\ref{fig:scalability_gaps}, the time required to find a high-quality, near-optimal solution (e.g., within a 5\% gap) is often a small fraction of the time needed to formally prove optimality. This allows managers to tailor the solver's runtime to the decision-making context, seeking guaranteed optimality for long-range planning while opting for rapid solutions for more tactical needs.

\subsubsection{Performance on Large-Scale Instances (N \texorpdfstring{$\geq$}{>=} 60)}
To define the practical limits of our approach, we tested the model on large-scale instances with a 3600-second time limit. Table~\ref{tab:large_scale_random_results} and Table~\ref{tab:large_scale_inc_results} show the results for the \texttt{RND} and \texttt{INC} instances, respectively, reporting the final cost and optimality gap for the MILP against the heuristic's performance.

\begin{figure}[htbp!]
    \centering
    \includegraphics[width=\textwidth]{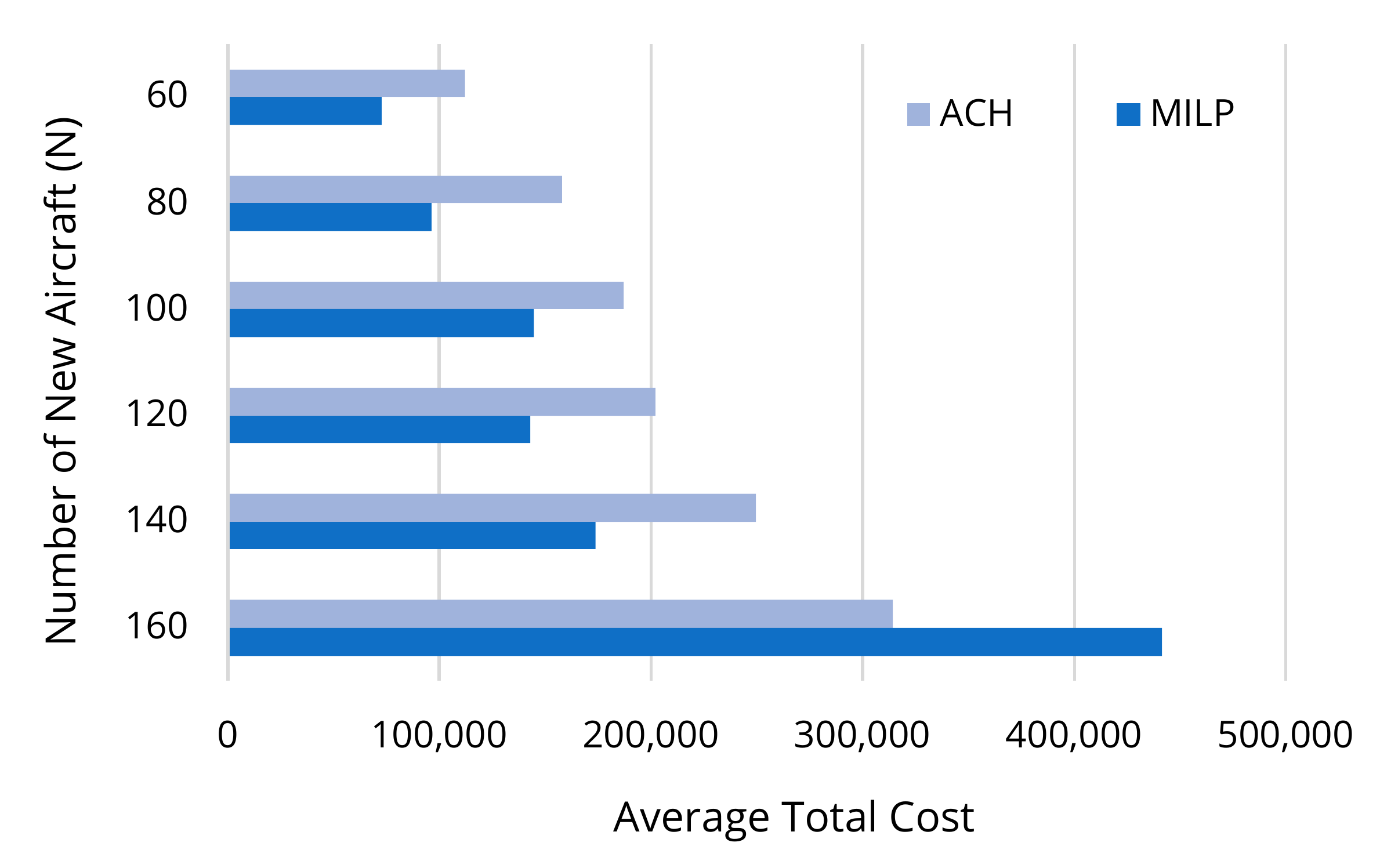} 
    \caption{Average total cost comparison for large-scale random (\texttt{RND}) instances under a 1-hour time limit. Note the performance inversion at N=160, where the heuristic finds a lower-cost solution.}
    \label{fig:large_scale_comparison}
\end{figure}

\begin{table}[htbp!]
\centering
\caption{Performance comparison on large-scale random instances (N $\geq$ 60) with a 3600-second time limit.}
\label{tab:large_scale_random_results}
\begin{tabularx}{\textwidth}{l *{6}{X}}
\toprule
\multirow{2}{*}{\thead{Instance Name}} & \multicolumn{2}{c}{\thead{MILP}} & \multicolumn{2}{c}{\thead{ACH Heuristic}} & \multicolumn{2}{c}{\thead{Accepted Aircraft}} \\
\cmidrule(lr){2-3} \cmidrule(lr){4-5} \cmidrule(lr){6-7}
& \thead{Total \\ Cost} & \thead{Final \\ Gap (\%)} & \thead{Total \\ Cost} & \thead{Time (s)} & \thead{ACH} & \thead{MILP} \\
\midrule
RND-N060-I01 & 77,721   & 4  & 131,267  & 63  & 27 & 39 \\
RND-N060-I02 & 86,267   & 6  & 113,555  & 66  & 25 & 36 \\
RND-N060-I03 & 53,737   & 18 & 91,298   & 67  & 27 & 40 \\
\midrule
RND-N080-I01 & 92,008   & 17 & 170,597  & 97  & 29 & 52 \\
RND-N080-I02 & 82,057   & 27 & 128,510  & 88  & 35 & 49 \\
RND-N080-I03 & 114,208  & 20 & 174,253  & 96  & 31 & 45 \\
\midrule
RND-N100-I01 & 154,695  & 44 & 176,116  & 120 & 37 & 43 \\
RND-N100-I02 & 154,016  & 39 & 191,847  & 114 & 41 & 37 \\
RND-N100-I03 & 125,545  & 38 & 192,985  & 119 & 37 & 55 \\
\midrule
RND-N120-I01 & 144,773  & 49 & 222,262  & 136 & 55 & 71 \\
RND-N120-I02 & 155,769  & 39 & 236,613  & 141 & 52 & 67 \\
RND-N120-I03 & 128,201  & 50 & 146,764  & 124 & 54 & 61 \\
\midrule
RND-N140-I01 & 147,768  & 45 & 218,253  & 159 & 60 & 76 \\
RND-N140-I02 & 168,740  & 58 & 247,066  & 170 & 59 & 70 \\
RND-N140-I03 & 205,143  & 47 & 283,302  & 171 & 57 & 73 \\
\midrule
RND-N160-I01 & 420,247  & 78 & 269,874  & 183 & 71 & 17 \\
RND-N160-I02 & 444,830  & 80 & 310,533  & 204 & 63 & 16 \\
RND-N160-I03 & 458,978  & 75 & 361,692  & 212 & 60 & 19 \\
\bottomrule
\end{tabularx}
\end{table}

\begin{table}[htbp!]
\centering
\caption{Performance comparison on large-scale incremental instances (N $\geq$ 60) with a 3600-second time limit.}
\label{tab:large_scale_inc_results}
\sisetup{group-separator={,}, group-minimum-digits=4}
\begin{tabular}{l S[table-format=6.0] S[table-format=2.0] S[table-format=6.0] S[table-format=3.0] S[table-format=2.0] S[table-format=2.0]}
\toprule
\multirow{2}{*}{\thead{Instance Name}} & \multicolumn{2}{c}{\thead{MILP}} & \multicolumn{2}{c}{\thead{ACH Heuristic}} & \multicolumn{2}{c}{\thead{Accepted Aircraft}} \\
\cmidrule(lr){2-3} \cmidrule(lr){4-5} \cmidrule(lr){6-7}
& {\thead{Total \\ Cost}} & {\thead{Final \\ Gap (\%)}} & {\thead{Total \\ Cost}} & {\thead{Time (s)}} & {\thead{ACH}} & {\thead{MILP}} \\
\midrule
INC-N060  & 85935  & 16 & 141040 & 83  & 19 & 35 \\
INC-N080  & 109527 & 28 & 174855 & 100 & 26 & 50 \\
INC-N100  & 126633 & 31 & 194227 & 122 & 35 & 58 \\
INC-N120  & 151017 & 31 & 229176 & 147 & 46 & 68 \\
INC-N140  & 198458 & 42 & 270716 & 176 & 51 & 77 \\
INC-N160  & 379240 & 70 & 297258 & 203 & 59 & 33 \\
\bottomrule
\end{tabular}
\end{table}

The analysis of larger instances tells a nuanced story. The \texttt{RND} instances reveal that problem size (N) is not the sole determinant of difficulty; the specific problem \textbf{structure} is paramount. This variance suggests that factors beyond sheer size, such as the temporal density of arrivals and departures or the geometric heterogeneity of the aircraft mix, define the true computational complexity. For most large-scale instances, the time-limited MILP continues to find solutions significantly better than the heuristic. However, at the largest scales (N=160), a notable inversion occurs where the heuristic produces a better solution, a phenomenon visualized in Figure~\ref{fig:large_scale_comparison}. This is a well-understood phenomenon in large-scale optimization: an exact solver, in its exhaustive search for a provable optimum, may not have sufficient time to navigate to a high-quality integer feasible solution within a vast search space. This finding provides a candid picture of the landscape: our MILP is the definitive tool for optimal planning up to very large scales, while a fast heuristic remains a viable fallback for extremely large, time-critical scenarios.

\subsection{Robustness to Temporal Congestion}
Finally, we return to the theme of congestion that motivated our work. Having seen our model's dramatic success on the congested benchmark case, this final experiment is designed as a controlled stress test to dissect the source of that success. To do this, we systematically increase the temporal density of aircraft arrivals for a fixed problem size (N=20). The results, shown in Table~\ref{tab:congestion_results}, quantify the model's performance degradation under this increasing pressure.

\begin{table}[htbp!]
\centering
\caption{Sensitivity analysis results for temporal congestion. All instances were solved to proven optimality (0\% gap).}
\label{tab:congestion_results}

\sisetup{group-separator={,}, group-minimum-digits=4}

\begin{tabular}{l S[table-format=3.0] S[table-format=5.0]}
\toprule
\thead{Instance Name} & {\thead{Total \\ Time (s)}} & {\thead{Total \\ Cost}} \\
\midrule
CON-N20-I01\_1.0x & 4   & 15231 \\
CON-N20-I01\_1.4x & 4   & 16396 \\
CON-N20-I01\_2.0x & 11  & 19905 \\
CON-N20-I01\_2.5x & 14  & 23709 \\
CON-N20-I01\_3.3x & 75  & 29597 \\
\midrule
CON-N20-I02\_1.0x & 2   & 16067 \\
CON-N20-I02\_1.4x & 11  & 26793 \\
CON-N20-I02\_2.0x & 14  & 30625 \\
CON-N20-I02\_2.5x & 90  & 33539 \\
CON-N20-I02\_3.3x & 100 & 38354 \\
\midrule
CON-N20-I03\_1.0x & 2   & 22222 \\
CON-N20-I03\_1.4x & 2   & 27949 \\
CON-N20-I03\_2.0x & 6   & 36303 \\
CON-N20-I03\_2.5x & 9   & 39117 \\
CON-N20-I03\_3.3x & 10  & 41962 \\
\bottomrule
\end{tabular}
\end{table}

The critical insight from Table~\ref{tab:congestion_results} is not that performance changes, but \textit{how} it changes. Unlike the brittle nature of event-based models, which can suffer catastrophic failure under high congestion, our model exhibits a \textit{graceful and controlled performance degradation}. While the solution time and cost logically increase with temporal density, the model continues to find optimal solutions in a computationally trivial timeframe. This test confirms that the fundamental design of our continuous-time model is inherently resilient to the very conditions that cripple alternative approaches, making it a robust and reliable tool for complex, real-world operational environments.

\section{Managerial Insights}
The results of our extensive computational study offer several critical insights for MRO management, translating the model's performance into actionable strategies for enhancing operational efficiency and financial returns.

\begin{itemize}
    \item \textbf{Driving Profitability Through Optimal Trade-off Management.}
    The model's primary value is its direct impact on financial performance. It achieves this by intelligently managing the crucial trade-off between maximizing hangar throughput (the number of serviced aircraft) and minimizing the systemic delay costs that additional aircraft can induce. This strategic optimization leads to transformative cost savings, with our results showing reductions of over \textbf{62\%} in total operational costs in representative cases (e.g., \texttt{RND-N010-I02}) compared to heuristic approaches. For managers, this capability to find the true economic optimum translates directly into enhanced profitability and a higher Return on Assets (ROA) from their critical hangar infrastructure.

    \item \textbf{Enhancing Strategic Agility with Near-Optimal Solutions.}
    For tactical planning, guaranteed optimality is not always required. Our results consistently show that solutions \textbf{within a 5\% optimality gap} are found in a fraction of the time needed to prove 0\% optimality, with negligible impact on cost. This provides managers with a crucial trade-off, enabling agile, day-to-day adjustments where speed is paramount, without sacrificing solution quality.

    \item \textbf{Quantifying the Financial Impact of Operational Congestion.}
    The model acts as a tool to price the cost of operational pressure. Our sensitivity analysis shows a direct correlation between temporal density and systemic cost; for example, a \textbf{3.3x increase in congestion} more than doubled the total costs for instance \texttt{CON-N20-I02}. Managers can use this capability for financial modeling to justify premiums for urgent requests and make more profitable slot allocation decisions.

    \item \textbf{Improving High-Stakes Strategic Decisions.}
    Our framework makes optimal solutions computationally accessible for complex problems previously deemed unsolvable. The ability to solve the congested benchmark \texttt{Case15-C9} in just \textbf{0.11 seconds} proves that managers no longer need to rely on heuristic answers for their most critical decisions. For financially significant planning, a guaranteed optimal plan mitigates risk and secures value that might otherwise be lost.

    \item \textbf{Choosing the Right Tool for Efficient Decision-Making.}
    The large-scale tests provide a nuanced understanding of the right tool for the job. While the MILP is superior for definitive planning, the result for instance \texttt{INC-N160}—where the heuristic found a solution approximately \textbf{21\% cheaper} than the time-limited MILP—highlights its value for rapid, large-scale exploratory analysis. The key for management is to establish a clear protocol for when to use each tool effectively.

\end{itemize}

These findings collectively demonstrate that advanced optimization is not just an operational tool, but a strategic enabler of revenue growth, cost control, and competitive advantage for MRO management.

\subsection{Solution Visualization}
To bridge the gap between the model's numerical output and practical application, a custom visualization tool was developed in Python. This tool serves as an interactive and dynamic decision support dashboard, translating the complex spatiotemporal solution into an intuitive graphical interface for hangar managers. Key features include a 2D animated layout of the hangar, an interactive timeline to step through events, and synchronized data tables.

\begin{sidewaysfigure}
    \centering
    \includegraphics[width=0.9\textheight]{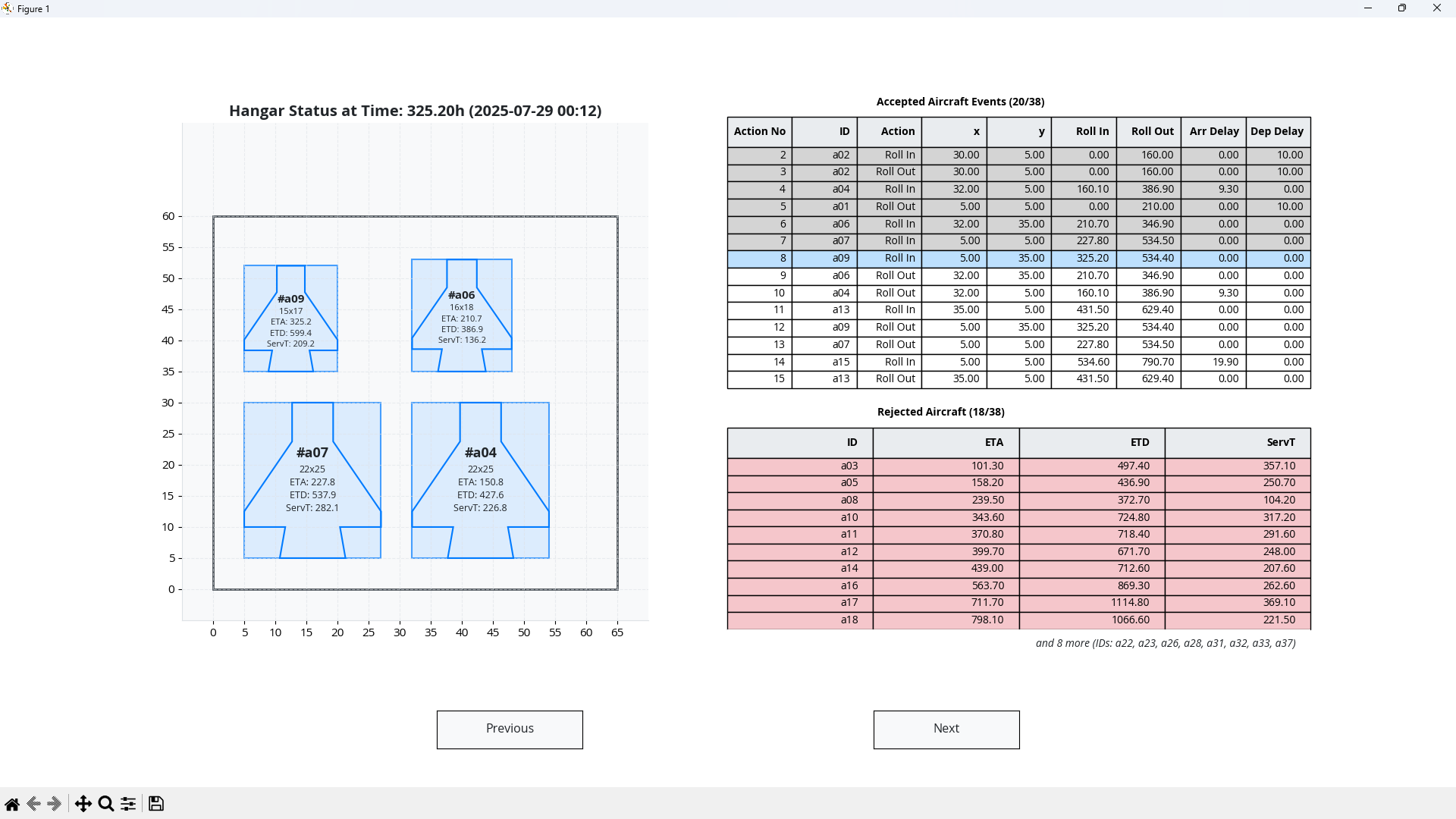} 
    \caption{A snapshot from the custom decision support and visualization tool, simulating a desktop application window. The interface displays the optimal hangar configuration for the INC-N036 test instance at time t=325.20h. The left panel shows the spatial layout of accepted aircraft, while the tables on the right provide detailed scheduling information for all accepted and rejected maintenance requests.}
    \label{fig:visualization}
\end{sidewaysfigure}

This tool is not merely for presentation; it is a powerful analytical instrument. As illustrated in Figure~\ref{fig:visualization}, it provides managers with a clear snapshot of the hangar's state at any point in time, enabling them to anticipate bottlenecks and understand the logic behind the optimal plan. The ability to visualize the intricate scheduling of movements—where the departure of one aircraft is precisely timed to clear a path for another—offers valuable managerial insights into the efficiency gains achievable through optimization. This transforms the abstract model into a tangible and actionable operational plan.

\section{Conclusion and Future Work}
This paper introduced a novel continuous-time MILP model that advances the state of the art in solving the integrated aircraft hangar scheduling and layout problem. Through an extensive and multi-faceted computational study, we demonstrated that this approach is not only methodologically rigorous but also practically powerful. Our model achieves orders-of-magnitude speedups on published benchmarks, solving a previously intractable congested instance in a fraction of a second. We established its scalability by finding provably optimal solutions for instances with up to 40 aircraft and providing verifiable bounds for problems as large as 160 aircraft. Furthermore, a controlled sensitivity analysis confirmed the model's robustness to temporal congestion, a critical bottleneck for prior event-based methods. For MRO managers, these results translate into actionable strategies for increasing hangar throughput, mitigating congestion costs, and supporting both tactical and strategic planning.

Future research could extend this model in several directions. A primary avenue is the incorporation of uncertainty. While our model provides a powerful deterministic backbone, integrating it with stochastic programming or robust optimization techniques would enhance its real-world applicability. Hybrid approaches that use our exact model for baseline planning, followed by real-time adjustment mechanisms, could offer both optimality and resilience. Such dynamic frameworks could draw inspiration from recent advances in reinforcement learning for long-term maintenance decisions \citep{hu2021reinforcement} or adaptive operating systems that manage the entire planning, scheduling, and execution (PSE) pipeline under uncertainty \citep{yang2025out}. Furthermore, our model's assumption of static placement could be relaxed to allow for multi-stage movements, addressing the complex blocking challenges modeled by \citet{chen2024stand}. Finally, enhancing the automated heuristic with a dynamic, look-back mechanism or other metaheuristic improvement phases could help close the performance gap while maintaining low computational times.

\section*{Declaration of Competing Interest}
The authors declare that they have no known competing financial interests or personal relationships that could have appeared to influence the work reported in this paper.

\section*{CRediT authorship contribution statement}
\textbf{Shayan Farhang Pazhooh:} Conceptualization, Methodology, Software, Data Curation, Formal analysis, Validation, Investigation, Writing – original draft, Writing – review \& editing, Visualization. \par

\textbf{Hossein Shams Shemirani:} Conceptualization, Methodology, Validation, Supervision, Project administration, Writing – review \& editing.

\section*{Declaration of Generative AI and AI-assisted Technologies in the Writing Process}
During the preparation of this work the author(s) used Gemini (a large language model by Google) in order to assist with language refinement and editing. After using this tool, the author(s) reviewed and edited the content as needed and take full responsibility for the content of the published article.

\section*{Funding}
This research did not receive any specific grant from funding agencies in the public, commercial, or not-for-profit sectors.

\section*{Data Availability Statement}
The data and source codes that support the findings of this study are openly available in Zenodo at \url{https://doi.org/10.5281/zenodo.17057809}. The development repository is also available on GitHub at \url{https://github.com/shayanfp/HangarModelResearch}.

\bibliographystyle{elsarticle-harv}
\bibliography{references}

\end{document}